\documentclass[a4paper,12pt]{article}

\usepackage{amsmath,amsfonts,amssymb}
\usepackage{amsthm}
\usepackage{indentfirst}
\usepackage{graphics}
\graphicspath{{Figures/}}

\theoremstyle{plain}
\newtheorem{lemma}{Lemma}[section]
\newtheorem{theorem}[lemma]{Theorem}

\theoremstyle{definition}
\newtheorem{definition}[lemma]{Definition}

\newtheorem{rem}[lemma]{Remark}

\DeclareMathOperator{\spec}{Spec}
\DeclareMathOperator{\spech}{\widehat{S}pec}
\DeclareMathOperator{\re}{Re}
\DeclareMathOperator{\im}{Im}
\DeclareMathOperator{\dom}{Dom}

\begin{document}

\title{Spectral pollution and second order relative spectra for
self-adjoint operators}

\author{Michael Levitin\\
\normalsize\small Department of Mathematics, Heriot-Watt University,\\
\normalsize\small Riccarton, Edinburgh EH14 4AS\\
\normalsize\small \texttt{M.Levitin@ma.hw.ac.uk}\\
\normalsize\small \texttt{www.ma.hw.ac.uk/$\sim$levitin}
\and
Eugene Shargorodsky\\
\normalsize\small School of Mathematical Sciences, University of Sussex,\\
\normalsize\small Falmer, Brighton, BN1 9QH, UK\\
\normalsize\small \texttt{E.Shargorodsky@sussex.ac.uk}\\
\normalsize\small \texttt{www.maths.sussex.ac.uk/Staff/ES/}}

\date{\today}

\bibliographystyle{amsplain}

\maketitle

\begin{abstract}
We consider the phenomenon of spectral pollution arising in calculation of spectra 
of self-adjoint operators by projection methods. We suggest a strategy of
dealing with spectral pollution by using the so-called second order relative spectra.
The effectiveness of the method is illustrated by a detailed analysis of two model 
examples.
\end{abstract}

{\small \textbf{Keywords:} computation of eigenvalues, spectral pollution, projection methods,
spurious eigenvalues, Agmon--Douglis--Nirenberg elliptic operators}

\

{\small \textbf{2000 Mathematics Subject Classification:} 47A75, 65J10}
\newpage

\tableofcontents

\listoffigures

\listoftables

\newpage

\section{Introduction}

Let $A$ be a self-adjoint operator  on a Hilbert space
$\mathcal{H}$. In order to determine the spectrum
$\spec(A)$ of
$A$, one can use the following projection method. Let
$(\mathcal{L}_k)_{k \in \mathbb{N}}$
be a sequence of closed linear subspaces of the domain
$\dom(A)$ of $A$ and suppose that the corresponding
orthogonal projections $P_k : \mathcal{H} \to \mathcal{L}_k$
converge strongly to the identity operator $I$. Let
$\Lambda(A)$ be the set of all such sequences of subspaces.
One may ask
whether or not
\begin{equation}\label{1}
\lim_{k \to \infty} \spec(A, \mathcal{L}_k) =
\spec(A) ,
\end{equation}
where  $\spec(A, \mathcal{L}_k)$ is the spectrum of $P_kA :
\mathcal{L}_k \to \mathcal{L}_k$ and ``$\lim$'' is defined in an
appropriate sense. Unfortunately the answer is, in general,
negative and $\lim\limits_{k \to \infty} \spec(A, \mathcal{L}_k)$
may contain points which do not belong to $\spec(A)$. This
phenomenon is well known in numerical analysis and is called {\em
spectral pollution}.

Since we allow $A$ to be unbounded, it is convenient for us to
extend $\spec(A)$ by adding to it $-\infty$, $+\infty$, or both,
if $A$ is unbounded below, above or from both sides, respectively.
We denote the resulting set by $\spech(A)$. In other words,
$\spech(A)$ is the closure of $\spec(A)$ in the two point
compactification $\overline{\mathbb{R}} :=
\mathbb{R}\cup\{\pm\infty\}$ of $\mathbb{R}$. Let
$\spec_{\text{ess}}(A)$ be the essential spectrum of $A$, i.e.
$$
\spec_{\text{ess}}(A) := \spec(A)\setminus\{\text{isolated
eigenvalues of finite multiplicity}\},
$$
and let $\spech_{\text{ess}}(A)$ be the
extended essential spectrum of $A$, i.e.
$$
\spech_{\text{ess}}(A) := \spech(A)\setminus\{\text{isolated
eigenvalues of finite multiplicity}\}.
$$

It turns out that spectral pollution may occur at any
point lying in a gap of the extended essential spectrum,
i.e. at any point of the set \
$\mbox{\rm conv}\left(\spech_{\text{ess}}(A)\right)\setminus
\spech_{\text{ess}}(A)$,
where ``conv'' means the convex hull (see the proof of Theorem
\ref{poll} below, where this fact is demonstrated
by choosing an appropriate
sequence of finite-dimensional subspaces).

One can argue that the
sequence $(\mathcal{L}_k)_{k \in \mathbb{N}}$ in the proof of
Theorem \ref{poll} had been deliberately chosen in a
particularly nasty way and that things like
this do not usually happen in ``real life''. This argument
is only partially true at best, and examples in Sections
3 and 4 show how the standard approach leads to
spectral pollution in some very straightforward
problems.

One can try to deal with spectral pollution by choosing
$\mathcal{L}_k$ in a ``reasonable way''. However the known
recipes (see, e.g., \cite{BBG}, \cite{BDG}, 
\cite{RSSV}, and references therein)
are by no means universal or rigorously justified for a wide
range of problems. An absolutely ``safe'' choice of
$(\mathcal{L}_k)_{k \in \mathbb{N}}$ would be the one based on
the spectral decomposition of $A$. Unfortunately knowing the
spectral decomposition is much more than knowing just the
spectrum which is to be found.

We propose a complementary (not an alternative!) way of dealing
with spectral pollution, which is based on the notion of the
second order spectrum $\spec_2(A, \mathcal{L})$ of $A$ relative to
a subspace $\mathcal{L} \subset \dom(A)$ (see Definition
\ref{Def}) introduced in \cite{EBD}. Second order relative spectra
usually contain complex numbers, but it turns out (see Theorem
\ref{mth}) that if $z \in \spec_2(A, \mathcal{L})$ then
\begin{equation}\label{eq:intersect}
\spec(A) \cap [\re z - |\im z|,\re z + |\im z|]\not= \varnothing.
\end{equation}
In particular, this means that if $z$ belongs to a second order
relative spectrum of $A$, and its imaginary part is small, than
there is a point of the spectrum of $A$ close to $\re z$.

\begin{figure}[thb!]
\begin{center}
\resizebox{0.9\textwidth}{!}{\includegraphics*{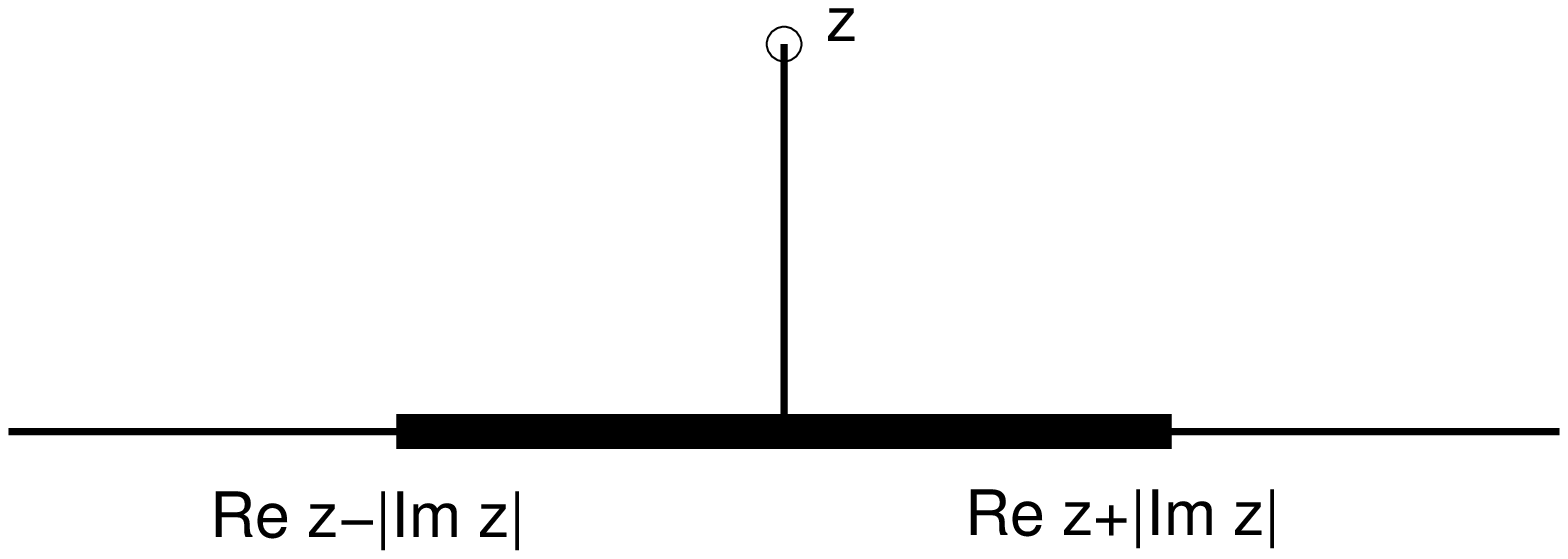}}
\end{center}
\caption{\label{fig:intersect}\small $\spec(A) \cap [\re z - |\im
z|,\re z + |\im z|]\not= \varnothing$.}
\end{figure}

This fact suggests the following strategy, which uses
\eqref{eq:intersect} as {\em a posteriori} estimate. {\em One
should choose a reasonable sequence $(\mathcal{L}_k)_{k \in
\mathbb{N}}$ and apply the standard projection method, i.e. find
$\spec(A, \mathcal{L}_k)$. Additionally, one should find
$\spec_2(A, \tilde{\mathcal{L}}_k)$ (for possibly different
$\tilde{\mathcal{L}}_k \subset \dom(A)$). Suppose $\spec(A,
\mathcal{L}_k)$ contain points which appear to converge to some
$\lambda \in \mathbb{R}$ as $k \to \infty$. If $\spec_2(A,
\tilde{\mathcal{L}}_k)$ also contain points lying close to
$\lambda$, then $\lambda$ is indeed close to the spectrum of $A$.
If however $\spec_2(A, \tilde{\mathcal{L}}_k)$ do not contain
points approximating $\lambda$, one cannot conclude that $\lambda$
is far from the spectrum. A definite claim which can be made in
this case is that $\lambda$ requires additional
attention\footnote{One can, for example, use the method suggested
in the companion paper by E.B. Davies and M. Plum \cite{DP}.}, as
spectral pollution may be around the corner.} It has to be said
though, that our (very limited!) experience seems to suggest that
$\spec_2(A, \tilde{\mathcal{L}}_k)$ do capture the spectrum quite
well.

We do not claim that the above method is going to be universally
successful. Our main ambition is merely to bring it to the attention
of those concerned with finding spectra numerically.
We have applied it to two problems which were chosen
mainly for their simplicity. The results (which turned out to be
better than we expected) are presented in Sections 3,
4, 6, and 7.

Section 2 contains the main theoretical results on standard (first
order) and second order relative spectra. Similar results on
second (and higher) order relative spectra have been obtained in
\cite{ES} under slightly stronger restrictions. We start with
Theorem  \ref{poll}, which has been mentioned above and which says
that spectral pollution may occur at any point lying in a gap of
the extended essential spectrum. It implies that one should be
particularly careful when applying a projection method to find the
spectrum of an operator $A$ in the following cases:\\
\qquad a) $A$ is not semi-bounded and $\spec(A) \not= \mathbb{R}$;\\
\qquad b) $A$ is semi-bounded but unbounded, and has a nonempty essential
spectrum different from a half-line;\\
\qquad c) $A$ is bounded and has at least one gap in the essential spectrum.\\
Examples of types
c) and b) are treated in Sections 3 and 4 respectively. We return
to these examples in  Sections 6 and 7 after giving the proofs of
the main theoretical results in Section 5.

Section 2 also discusses an
additional difficulty which arises for unbounded operators.
It is well known that if $A$ is bounded, then
the left-hand side of \eqref{1} contains the right-hand side.
If, however, $A$ is unbounded above, say, then
the left-hand side of \eqref{1}
may well be equal to $\{+\infty\}$ (see Theorem \ref{unb}).

In order to make the paper more readable, we have divided
the material into two parts. The first one consists of
Sections 2--4 and covers the main theoretical and numerical
results. The same material is treated at a greater depth
in the second part of the paper, which consists of Sections
5--7. The summary is given in Section 8.

Finally, the accompanying webpage\\
\centerline{\texttt{www.ma.hw.ac.uk/$\sim$levitin/SpectralPollution}}
contains the listings
of the Matlab programmes written by the authors and used for numerical 
calculations in Examples I and II.

\section{Main results}

Let $A$ be a self-adjoint operator  on a Hilbert space
$\mathcal{H}$. Throughout the paper $E(\cdot)$ denote the spectral
projections associated with $A$. The proofs of all results of
this Section are given in Section 5.

The following result states that spectral pollution may occur at {\em
any} point in a gap of the extended essential spectrum of $A$.

\begin{theorem}\label{poll}
For any $ \lambda \in
\textup{conv}\left(\spech_{\textup{ess}}(A)\right)\setminus
\spech_{\textup{ess}}(A)$
there exists an increasing sequence
$(\mathcal{L}_k)_{k \in \mathbb{N}} \in \Lambda(A)$ such that
$$
\lambda \in \spec(A, \mathcal{L}_k) , \qquad \forall k
\in\mathbb{N} .
$$
\end{theorem}

It is well known that for bounded self-adjoint operators
the left-hand side of \eqref{1} contains the right-hand side:
\begin{equation}\label{upp}
\lim_{k \to \infty} \spec(A, \mathcal{L}_k) \supseteq
\spec(A)
\end{equation}
(see, e.g.,  \cite{Arv} or \cite{ES}). So, spectral pollution
means that the left-hand side of \eqref{1} is strictly larger
than the right-hand side. For unbounded self-adjoint operators
the situation is more complicated: the left-hand side of \eqref{1}
does not necessarily contain the right-hand side. The reason
for this is intuitively very clear. Let $\mathcal{L}_k$ be
finite dimensional and suppose for definiteness that
$A$ is not bounded above. Then perturbing a basis of $\mathcal{L}_k$
by arbitrarily small vectors lying in
$\dom(A)\cap E([N, +\infty))\mathcal{H}$ with sufficiently large $N$,
one can push the spectrum of the corresponding finite dimensional
operator arbitrarily far to the right.

\begin{theorem}\label{unb}
Let $A$ be an unbounded above self-adjoint operator on a separable
Hilbert space $\mathcal{H}$. Then for an arbitrary sequence
$(\mathcal{L}_k)_{k \in \mathbb{N}} \in \Lambda(A)$ of finite
dimensional subspaces and arbitrary sequences
$(\varepsilon_k)_{k\in \mathbb{N}} \searrow 0$ and
$(R_k)_{k\in\mathbb{N}} \nearrow +\infty$, there exists a sequence
$(\mathcal{L}'_k)_{k \in \mathbb{N}} \in \Lambda(A)$ such that
$$
\|P_k - P'_k\| < \varepsilon_k \qquad\text{and}\qquad\spec(A,
\mathcal{L}'_k) \subset (R_k, +\infty) , \quad\forall k \in
\mathbb{N} ,
$$
where $P'_k : \mathcal{H} \to \mathcal{L}'_k$ are the corresponding
orthogonal projections. A similar statement holds for operators
unbounded below.
\end{theorem}

The requirement $\mathcal{L}_k \subset \dom(A)$
is too restrictive for many applications, in particular
for approximating spectra of differential operators by
finite element methods. When using finite element methods,
one usually states the problem in a weak (variational)
form and requires that $\mathcal{L}_k$ belong to the
domain of the corresponding sesquilinear form. The
restriction $\mathcal{L}_k \subset \dom(A)$
would mean extra smoothness of the finite elements.

Suppose $A$ is a bounded below
self-adjoint operator on a separable Hilbert space $\mathcal{H}$
and let $Q_A : \mathcal{H}\times\mathcal{H} \to
\mathbb{C}$ be the corresponding closed sesquilinear form
(see, e.g., \cite[\S 1, Ch. VI]{Kato}).
Let
$(\mathcal{L}_k)_{k \in \mathbb{N}}$
be a sequence of finite dimensional linear subspaces of the domain
$\dom(Q_A)$ of $Q_A$ such that the corresponding
orthogonal projections $P_k : \mathcal{H} \to \mathcal{L}_k$
converge strongly to $I$. Let
$\Lambda(Q_A)$ be the set of all such sequences.
A number $\lambda \in \mathbb{C}$ is said
to belong to $\spec(A, \mathcal{L}_k)$
if there exists $u \in \mathcal{L}_k\setminus\{0\}$ such that
$$
Q_A(u, v) - \lambda (u, v) = 0 , \ \
\forall v \in \mathcal{L}_k .
$$
If $\mathcal{L}_k \subset \dom(A)$, then $Q_A(u, v) = (Au, v)$ and
the above definition means that\\
$\bullet$ there exists $u \in
\mathcal{L}_k\setminus\{0\}$ such that $(Au, v) - \lambda (u, v) =
0$, $\forall v \in \mathcal{L}_k$,\\
$\bullet$ i.e. there exists
$u \in \mathcal{L}_k\setminus\{0\}$ such that $P_kAu = \lambda u$\\ 
$\bullet$ i.e. $\lambda$ belongs to the spectrum of $P_kA :
\mathcal{L}_k \to \mathcal{L}_k$.\\
Hence one gets the old
definition of $\spec(A, \mathcal{L}_k)$ for finite dimensional
$\mathcal{L}_k \subset \dom(A)$ (cf. Remarks~\ref{spec2}
and~\ref{rem:saidbefore}(b) below).

It is not difficult to prove that Theorem~\ref{unb} remains true
if $\Lambda(A)$ is replaced by $\Lambda(Q_A)$ (see Section 5).

The above results show the difficulties one may
face when trying to approximate the spectrum of
a self-adjoint operator with the help of a
projection method. A possible way of addressing
spectral pollution is based on the notion of
the second order relative spectrum.

\begin{definition}\label{Def}
Let $A$ be a self-adjoint operator on a Hilbert space $\mathcal{H}$
and let $\mathcal{L}$ be a finite dimensional subspace of
$\dom(A)$. A number $z \in \mathbb{C}$ is said
to belong to the {\em second order spectrum
$\spec_2(A, \mathcal{L})$ of $A$ relative to $\mathcal{L}$}
if there exists $u \in \mathcal{L}\setminus\{0\}$ such that
\begin{equation}\label{def}
\left((A - z I)u, (A - \overline{z} I)v\right) = 0 , \ \
\forall v \in \mathcal{L} ,
\end{equation}
where $(\cdot, \cdot)$ is the scalar
product in $\mathcal{H}$.
\end{definition}

\begin{rem}\label{spec2}
Suppose $\mathcal{L} \subset \dom(A^2)$
and let $P$ be the orthogonal projection onto $\mathcal{L}$.
Then $z \in \mbox{\rm Spec}_2(A, \mathcal{L})$ \\
$\bullet$ iff there exists
$u \in \mathcal{L}\setminus\{0\}$ such that
$\left((A - z I)^2 u, v\right) = 0$,
$\forall v \in \mathcal{L}$, \\
$\bullet$ i.e. iff there exists
$u \in \mathcal{L}\setminus\{0\}$ such that
$u \in \mbox{\rm Ker}(P(A - z I)^2)$, \\
$\bullet$ i.e. iff the operator
$P(A - z I)^2 : \mathcal{L} \to \mathcal{L}$
is not invertible. \\
The last is the definition of the second order relative spectrum
introduced in \cite{EBD} and used in \cite{ES}.
\end{rem}

\begin{theorem}\label{mth}
If $z \in \spec_2(A, \mathcal{L})$ then
$$
\spec(A) \cap [\re z - |\im z|,
\re z + |\im z|]
\not= \varnothing .
$$
\end{theorem}

\begin{rem}\label{rem:saidbefore}
(a) Theorem \ref{mth} has been proved in
\cite{ES} in the case
$\mathcal{L} \subset \dom(A^2)$. \\
(b) Everything said above about second order relative spectra
remains unchanged if $\mathcal{L}$
is infinite dimensional, with only one exception: the implication
$$
P(A - z I)^2 : \mathcal{L} \to \mathcal{L} \text{ is not
invertible }\ \Longrightarrow \ \ \text{Ker}(P(A - z I)^2) \cap
\mathcal{L} \not= \{0\}
$$
(see Remark \ref{spec2}) is no longer true.
\end{rem}

\section{Example I: operator of multiplication
by a discontinuous function}

Consider the operator of multiplication by a
function $a$:
$$
A = aI : L_2([-\pi, \pi]) \to L_2([-\pi, \pi])\,,
$$
where
$$
a(x) =
\begin{cases}
\displaystyle -\frac{3}{2}+\frac{1}{2}\cos\sqrt{5}\,x\,,
&\qquad\text{for}\quad -\pi\le x< 0\,,\\
\displaystyle 2+\cos\sqrt{2}\,x\,,
&\qquad\text{for}\quad 0\le x<\pi\,.
\end{cases}
$$
It is clear that
$$
\spec(A) = \spec_{\text{ess}}(A) =
a\left([-\pi, \pi]\right) = [-2, -1] \cup [1, 3]\, .
$$
The function $a$ and, in
particular, the irrational factors $\sqrt{5}$ and $\sqrt{2}$
have been chosen to get a ``generic" situation and avoid
any possible hidden symmetries.

\begin{figure}[thb!]
\begin{center}
\resizebox{0.9\textwidth}{!}{\includegraphics*{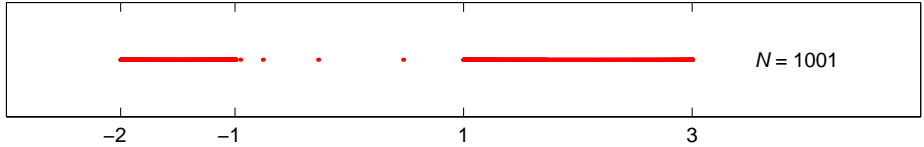}}
\end{center}
\caption{\label{fig:spec_A_L_N}\small $\spec(A,\mathcal{L}_{1001})$.}
\end{figure}

Using the standard orthonormal basis $e_k(x) := (2\pi)^{-1/2}
\exp{(ikx)}$, $k \in \mathbb{Z}$, one can represent $A$ as a
(doubly) infinite Toeplitz matrix and take the orthogonal
projections onto $\mathcal{L}_N := \text{span}\,(\{e_k\}_{k =
-n}^n)$, where $N=\text{dim}\,\mathcal{L}_N=2n+1$.
Figure~\ref{fig:spec_A_L_N} shows the approximation of the
spectrum obtained by the standard projection method for a typical
value of $N$. One can see that the essential spectrum $[-2, -1]
\cup [1, 3]$ is approximated very well, but that, additionally,
there are some computed eigenvalues lying in the gap of the
essential spectrum.

The spurious eigenvalues vary very slowly with $N$, see
Table~\ref{tab:spec_A}, and their number
$\#(\spec(A,\mathcal{L}_N)\cap(-1,1))$ also increases slowly as
$N$ increases. In fact, they seem to converge monotonically to
some points in $(-1,1)$, and judging by the numerical evidence
alone, one might have thought that the limits belonged to the
spectrum, if the example was less trivial and the spectrum was not
known.

\begin{table}[thb!]
\begin{center}
\begin{tabular}[t]{|c|*{4}{r}|}
\hline $N$ & 101  & 401  & 701   & 1001\\ \hline
 & & &  -0.9997 & -0.9985\\
         &  &  -0.9834 & -0.9673 & -0.9550 \\
$\spec(A,\mathcal{L}_N)\cap(-1,1)$ &  -0.9264&  -0.8250 & -0.7823
& -0.7553 \\
    & -0.4866 & -0.3478 & -0.3002& -0.2721\\
    &  0.4362 &  0.4597  & 0.4673 & 0.4717\\ \hline
\end{tabular}
\end{center}
\caption{\label{tab:spec_A}\small Eigenvalues lying in the gap of
the essential spectrum of $A$ for some values of $N$.}
\end{table}

The presence of spurious points is hardly surprising. Indeed, if
instead of the operator of multiplication by $a$, one considers
the corresponding Toeplitz operator $T_a$ acting on the Hardy
space $\mathtt{H}^2 =\mbox{span}\,(\{e_k\}_{k \ge 0})$, then the projection 
method using $\mathcal{M}_N := \mbox{span}\,(\{e_k\}_{k = 0}^{N - 1})$ 
leads to the \emph{same}
finite-dimensional Toeplitz matrices. In other words, the standard
projection method does not know which spectrum is being
calculated: that of  the operator of multiplication or of the
Toeplitz operator. The latter equals $[-2, 3]$ (see \cite{HW} or
\cite[2.36]{BS}). Therefore not only will we have spectral
pollution, the spurious eigenvalues will eventually fill the
entire gap $(-1, 1)$ (cf. \eqref{upp} or \cite{Boet}).

\begin{rem}\label{slow}
It is interesting that spurious eigenvalues appear so slowly.
Three spurious eigenvalues are present already for rather small
values of $N$. The fourth one emerges from $-1$ for $N \approx
161$, and the fifth one emerges from $-1$ for $N \approx 631$.
Such a slow rate seems to have the following explanation. The part
$(-1, 1)$ of the spectrum of the Toeplitz operator $T_a$ is there
because of the discontinuities of $a(x)$ at $x=0$ and $x=-\pi \sim
\pi$. Functions responsible for this part of the spectrum appear
to be very localized at $0$ and $-\pi \sim \pi$ and hence have a
very slow decay of the Fourier coefficients.
Figure~\ref{fig:typical_ef} shows the graph of the trigonometric
polynomial which is an eigenfunction of $P_N A : \mathcal{L}_N \to
\mathcal{L}_N$, $N = 1001$, corresponding to $\lambda\approx
0.4717$.

\begin{figure}[thb!]
\begin{center}
\resizebox{0.9\textwidth}{!}{\includegraphics*{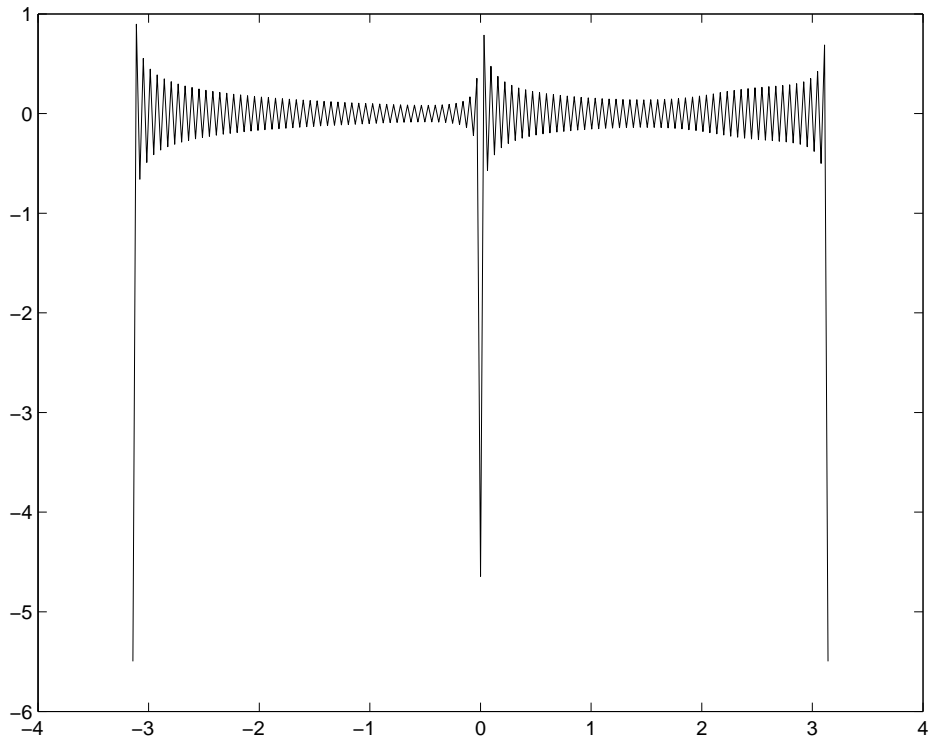}}
\end{center}
\caption{\label{fig:typical_ef}\small A typical computed
``eigenfunction'' of $A$ corresponding to a spurious eigenvalue.}
\end{figure}
\end{rem}

We now show how second order relative spectra
locate the spectrum of $A$ and, in particular,
detect spectral pollution in the standard
projection method.

\begin{figure}[thb!]
\begin{center}
\resizebox{0.9\textwidth}{!}{\includegraphics*{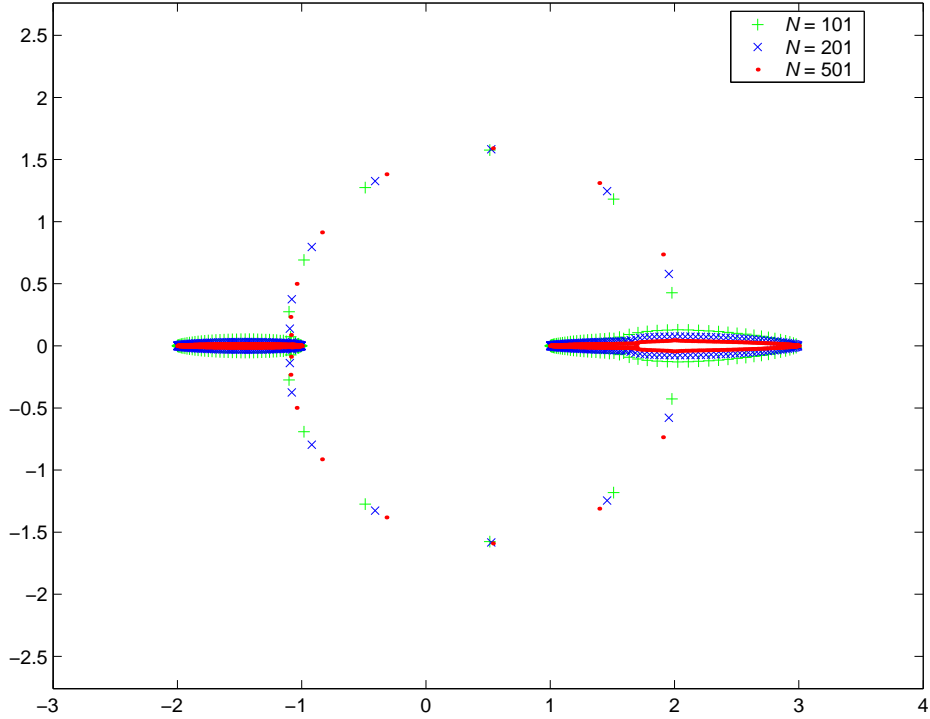}}
\end{center}
\caption{\label{fig:spec2_A}\small $\spec_2(A, \mathcal{L}_N)$.}
\end{figure}

Figure~\ref{fig:spec2_A} shows the second order relative spectra
$\spec_2(A, \mathcal{L}_N)$ computed for $N = 101, 201, 501$.

For each $N$, the second order relative spectrum consists
of a few points
lying approximately on a circle centred on the real axis
and two contours ``surrounding'' two intervals of the (essential)
spectrum. The number of points of $\spec_2(A, \mathcal{L}_N)$
lying on the circle does not seem to increase with $N$.

\begin{figure}[thb!]
\begin{center}
\resizebox{0.9\textwidth}{!}{\includegraphics*{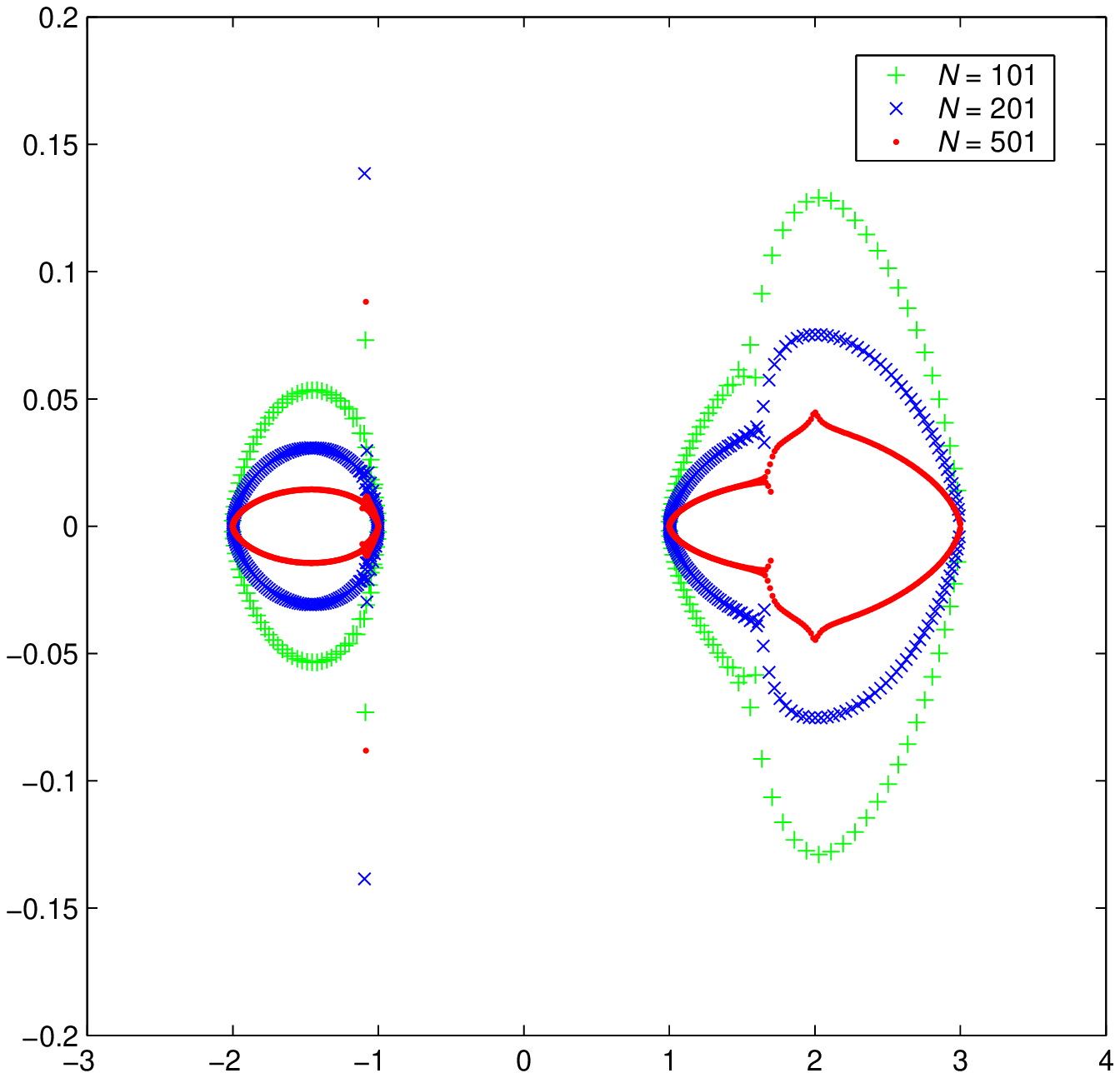}}
\end{center}
\caption{\label{fig:spec2_A_detail}\small Detailed view of
Figure~\ref{fig:spec2_A} near the real axis.}
\end{figure}

Figure~\ref{fig:spec2_A_detail} gives a detailed view of
Figure~\ref{fig:spec2_A} near the real axis (note the different
scales along the coordinate axes). The contours made of the points
of the second order relative spectra, which ``surround'' two
intervals of the (essential) spectrum, get closer to the real line
as $N$ increases. Most importantly, the absence of points of the
second order relative spectra close to the eigenvalues found by
the standard projection method inside the interval $(-1,1)$
suggests that these eigenvalues are spurious.

In order to illustrate a situation when an operator has both
non-empty essential and discrete spectra, let us consider a
two-dimensional perturbation of $A = aI$,

\begin{equation}\label{eq:B}
\begin{split}
Bu(x) := a(x)u(x) &- \frac{3\exp{(-ix)}}{2\pi} \int_{-\pi}^\pi
u(t)\exp{(it)} dt \\ &+ \frac{\exp{(2ix)}}{2\pi} \int_{-\pi}^\pi
u(t)\exp{(-2it)} dt \, , \qquad x \in [-\pi, \pi] .
\end{split}
\end{equation}

\begin{figure}[thb!]
\begin{center}
\resizebox{0.9\textwidth}{!}{\includegraphics*{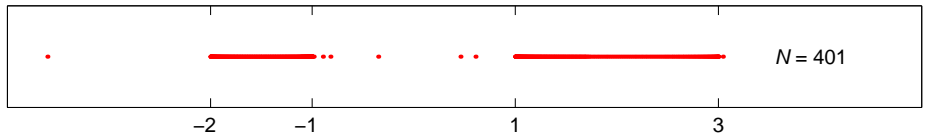}}
\end{center}
\caption{\label{fig:spec_B_L_N}\small $\spec(B,\mathcal{L}_{401})$.}
\end{figure}

This time, the essential spectrum is again very well approximated
by the usual projection method, and the perturbation gives rise to
several eigenvalues, see Figure~\ref{fig:spec_B_L_N}, listed here
for $N=401$ : one ($\lambda_1\approx -3.6057$) to the left of the
essential spectrum, six ($\lambda_2\approx -0.9842$,
$\lambda_3\approx -0.8947$, $\lambda_4\approx -0.8203$,
$\lambda_5\approx -0.3496$, $\lambda_6\approx 0.4603$, and
$\lambda_7\approx 0.6082$) in the gap $(-1, 1)$, and one
$\lambda_8\approx 3.0432$) to the right of the essential spectrum.

\begin{figure}[thb!]
\begin{center}
\resizebox{0.9\textwidth}{!}{\includegraphics*{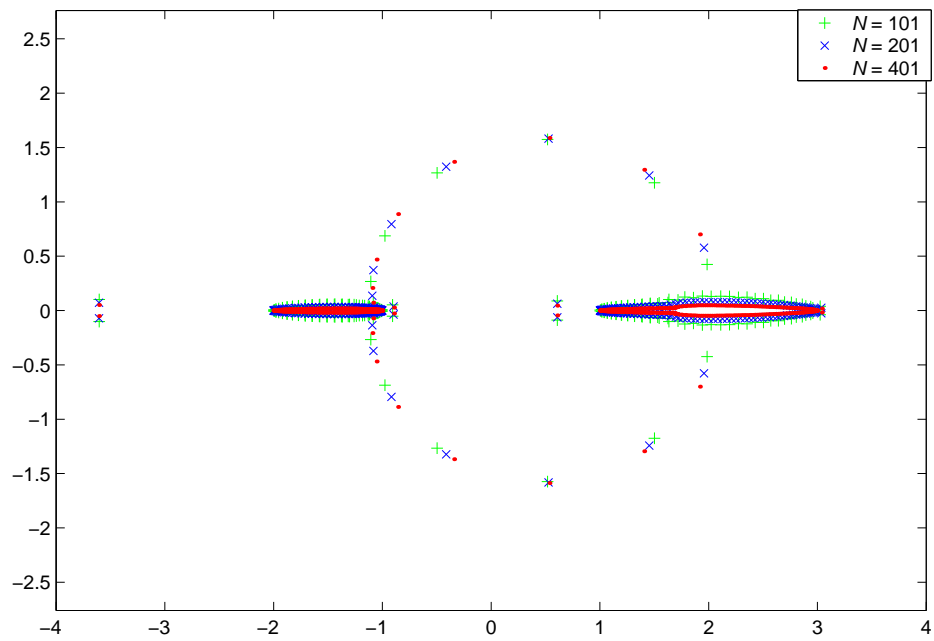}}
\end{center}
\caption{\label{fig:spec2_B}\small $\spec_2(B, \mathcal{L}_N)$.}
\end{figure}

\begin{figure}[thb!]
\begin{center}
\resizebox{0.9\textwidth}{!}{\includegraphics*{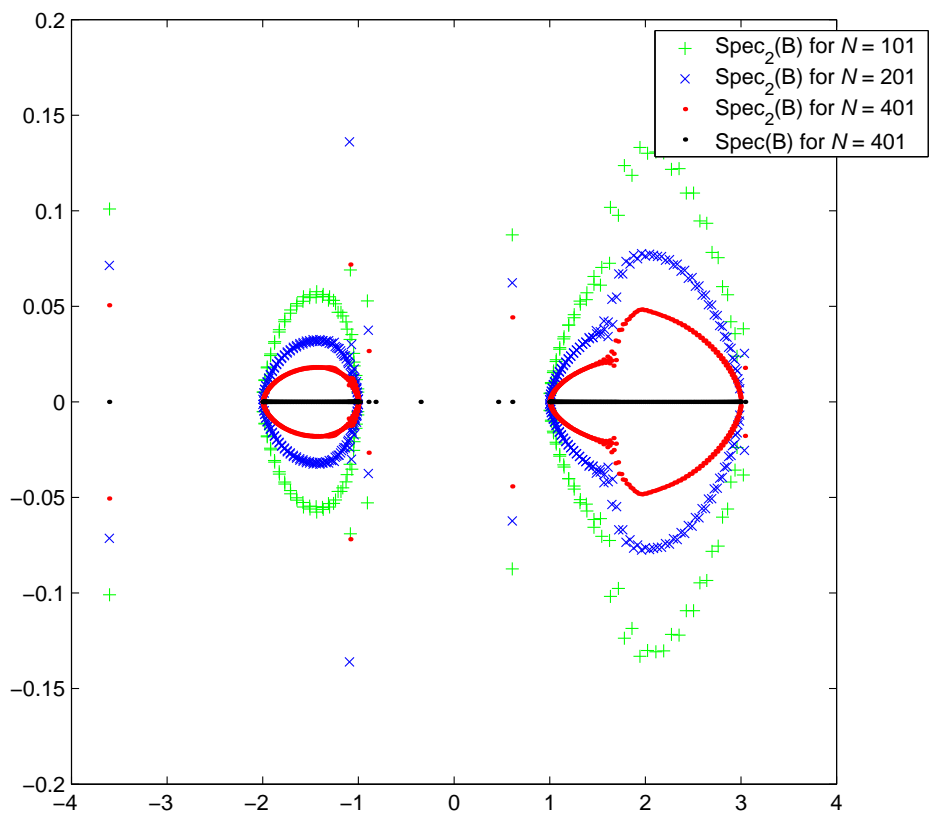}}
\end{center}
\caption{\label{fig:spec2_B_detail}\small Detailed view of
Figure~\ref{fig:spec2_B} near the real axis and, for comparison,
$\spec(B, \mathcal{L}_N)$.}
\end{figure}

To establish which of the eight eigenvalues are spurious, we
compute the second order relative spectra of $B$. Their structure
is similar to that of $\spec_2(A, \mathcal{L}_N)$ apart from
several additional points located very close to the real axis, see
Figures \ref{fig:spec2_B} and \ref{fig:spec2_B_detail}. In
particular, for $N=401$ these points are approximately $-3.6056
\pm 0.0505i$, $-0.8929 \pm 0.0266i$, $0.6085 \pm 0.0442i$, and
$3.0413 \pm 0.0178i$. Thus, the eigenvalues $\lambda_1$,
$\lambda_3$, $\lambda_7$ and $\lambda_8$ are correct, whereas the
other eigenvalues $\lambda_2$, $\lambda_4$, $\lambda_5$, and
$\lambda_6$ are almost certainly spurious.

\clearpage

\section{Example II: Stokes type system of ODEs}

Our next example is an operator $H$ in $L_2([-\pi,\pi])$ generated
by a simple $2\times2$ Stokes type system of ODEs subject to the
periodic boundary conditions,
$$
H=
\begin{pmatrix}
-\dfrac{d^2}{dx^2} & -\dfrac{d}{dx}\\[+1em]
 \dfrac{d}{dx} & 2
\end{pmatrix}\, .
$$
Finding the spectrum of $H$ explicitly is easy. The essential
spectrum is found by using a general result on
Agmon-Douglis-Nirenberg (that is, mixed order) elliptic systems
\cite{GG}, which gives
$$
\spec_{\text{ess}}(H) = \left\{\lambda \in \mathbb{C} \Big| \
\exists \xi \in \mathbb{R}\setminus\{0\} : \
\det \begin{pmatrix}
\xi^2  & i\xi\\
 -i\xi & 2 - \lambda
\end{pmatrix} = 0\right\}
= \{1\} .
$$
To find the discrete spectrum, consider the system of
equations
\begin{equation}\label{np1}
\left\{\begin{alignedat}{3}
-u''&-v'&&=\lambda u\,,\\
u'&+2v&&=\lambda v
\end{alignedat}\right.
\end{equation}
subject to the periodic boundary conditions. Expressing $v$ from
the second equation and substituting the result into the first one
gives
\begin{equation}\label{eq:sixstar}
-u''=\frac{\lambda(\lambda -2)}{\lambda-1}\,u\,.
\end{equation}

Comparing this equation with the standard periodic spectral
problem for $-d^2/dx^2$ we conclude that
\begin{equation}\label{eq:eigsH}
\lambda\in\spec_{\text{disc}}(H)\quad\text{iff}\quad
\frac{\lambda(\lambda -2)}{\lambda-1}\in
\spec\left(-\frac{d^2}{dx^2}\right)=\{k^2: k\in\mathbb{Z}\}\,.
\end{equation}

This implies that $H$ has
two series of eigenvalues $\lambda_k^\pm$ given by
$$
\lambda_k^\pm=\frac{k^2+2\pm\sqrt{(k^2+2)^2-4k^2}}{2}\,,\qquad
k=0,1,2,\dots
$$
Each eigenvalue $\lambda_k^-$ or $\lambda_k^+$ with $k>0$ is of
multiplicity 2. The eigenvalues of the ``minus'' series are
located in the interval $[0,1)$ and converge to $1$ as
$k\to\infty$. The eigenvalues of the ``plus'' series are located
in the interval $[2,+\infty)$ and tend to $+\infty$ as
$k\to\infty$. Note that the interval $(1,2)$ is spectrum-free. The
eigenvectors $(u,v)$ of $H$ corresponding to the eigenvalue
$\lambda_k^\pm$, $k>0$ are linear combinations of
\begin{equation}\label{eigv}
\begin{pmatrix}
1\\
ik/(\lambda_k^\pm - 2)\end{pmatrix} e^{ik x}
\ \ \text{and} \ \
\begin{pmatrix}
1\\
-ik/(\lambda_k^\pm - 2)\end{pmatrix} e^{-ik x} .
\end{equation}
The eigenvectors corresponding to $\lambda_0^- = 0$
and $\lambda_0^+ = 2$ are proportional to
$\begin{pmatrix}
1\\
0\end{pmatrix}$ and $\begin{pmatrix}
0\\
1\end{pmatrix}$ respectively.

The projection method we use in the numerical analysis of the
spectrum of $H$ is the standard finite element method (FEM). It
turns out (see Section 7 below) that the results of calculations
depend dramatically upon the choice of the corresponding spaces of
finite elements.

We start with the simplest possible choice --- the space
$\mathbf{L}^p_N$ of one dimensional periodic Lagrange elements of
order $p$, that is, of continuous functions which are piecewise
polynomials of degree $p$. Here $N$ is the number of mesh node
points on the interval $[-\pi,\pi]$, so the total dimension of
$\mathbf{L}^p_N$ is $Np$.

\begin{rem} In most of the numerical experiments, we use just the linear
Lagrange elements (i.e., of order $1$). We shall therefore write
$\mathbf{L}_N:=\mathbf{L}^1_N$. Additional tests have shown that
varying the order affects neither the qualitative picture, nor the
rate of convergence.

It is important to note that the spaces $\mathbf{L}^p_N$ are
embedded into the Sobolev space $W_2^1([-\pi,\pi])$ but not into the
Sobolev spaces $W_2^2([-\pi,\pi])$. Thus, the elements of
$\mathbf{L}_N\times\mathbf{L}_N$ do not necessarily belong to
$\dom(H)$. However, they belong to the domain of the corresponding
quadratic form, which is sufficient for the FEM.
\end{rem}

The Lagrange elements  are provided by all standard FEM packages,
including, e.g., FEMLAB (see \cite{Femlab}).

\begin{figure}[thb!]
\begin{center}
\resizebox{0.9\textwidth}{!}{\includegraphics*{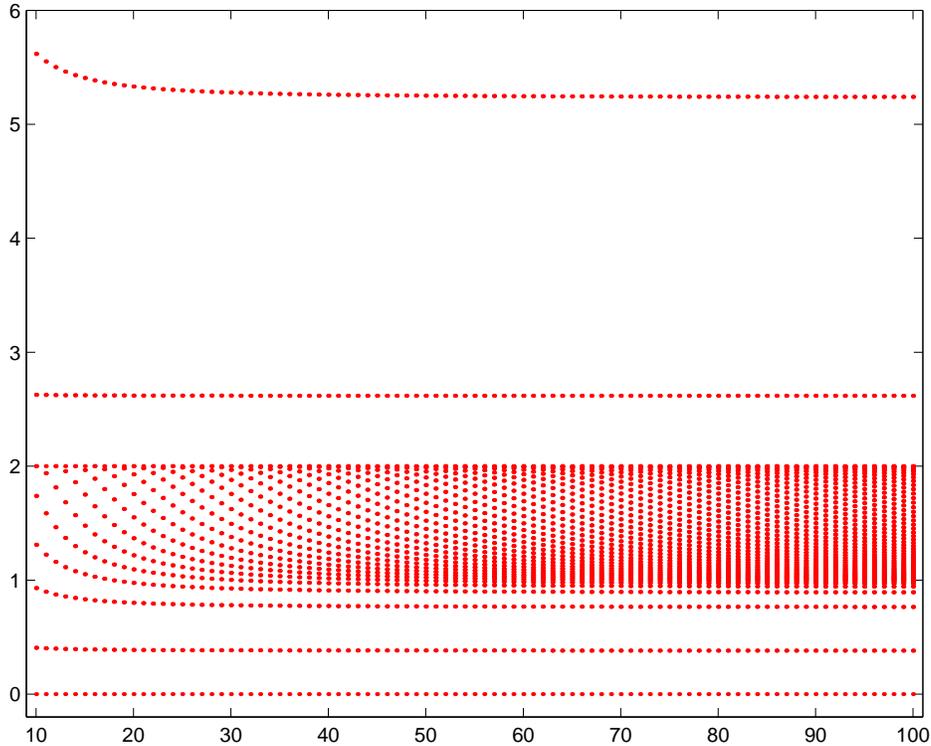}}
\end{center}
\caption{\label{fig:spec_H}\small FEM calculated
$\spec(H,\mathbf{L}_N\times\mathbf{L}_N)$ plotted against $N$. All
eigenvalues in the interval $(1,2)$ are spurious!}
\end{figure}

The calculation of $\spec(H,\mathbf{L}_N\times\mathbf{L}_N)$
produces a very striking result: it gives very accurately several
initial eigenvalues from the $\lambda^\pm$ series, filling at the
same time the entire interval $(1,2)$ with spurious eigenvalues
(see Figure~\ref{fig:spec_H}, which, for brevity, shows only the
eigenvalues lying in the interval $[0,6]$). In the particular case
$N=100$, there are 15 eigenvalues in $[0,1]$, 99 eigenvalues in
$[2,\infty)$, and the remaining 86 eigenvalues are in $(1,2)$.

The subspace $\mathcal{L}=\mathbf{L}_N\times\mathbf{L}_N$ used
above is unsuitable for  the calculation of the second order
relative spectra $\spec_2(H, \mathcal{L})$ since for elements of
this space the left-hand side of \eqref{def} does not make sense.
Therefore, as usual when dealing with fourth order ODES, let us
introduce the standard space of periodic cubic Hermite elements
$\mathbf{M}_N$, whose functions are {\em continuously
differentiable\/} and are piecewise polynomials of degree 3. They
belong to the Sobolev space $W_2^2([-\pi,\pi])$ so they are
smoother than elements of $\mathbf{L}_N$. Here again $N$ is the
number of mesh node points, so that $\dim \mathbf{M}_N = 2N$.

\begin{figure}[thb!]
\begin{center}
\resizebox{0.9\textwidth}{!}{\includegraphics*{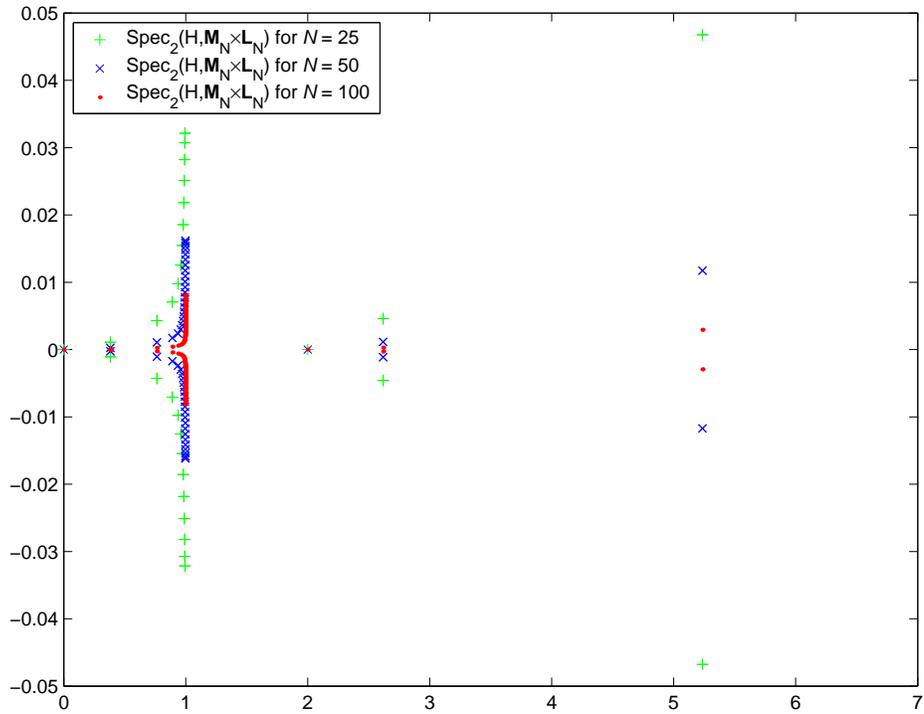}}
\end{center}
\caption{\label{fig:spec2_H_detail}\small Second order relative
spectra $\spec_2(H, \mathbf{M}_N\times\mathbf{L}_N)$ for various
values of $N$ --- low eigenvalues.}
\end{figure}

\begin{figure}[thb!]
\begin{center}
\resizebox{0.9\textwidth}{!}{\includegraphics*{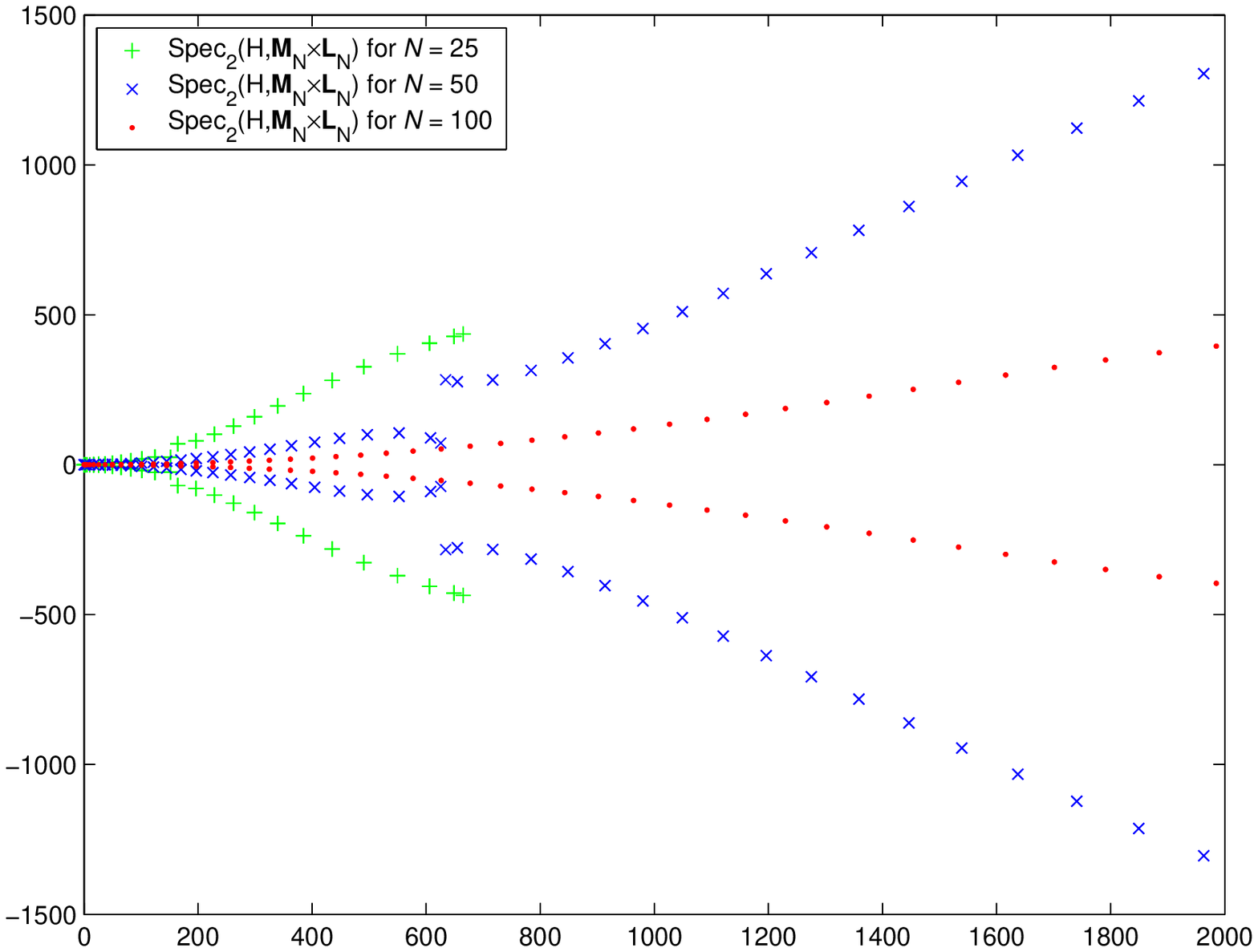}}
\end{center}
\caption{\label{fig:spec2_H}\small Second order relative spectra
$\spec_2(H, \mathbf{M}_N\times\mathbf{L}_N)$ for various values
of $N$ --- all eigenvalues.}
\end{figure}

Since the operator $H$ is of mixed order, it is sufficient to
choose the space $\mathbf{M}_N\times\mathbf{L}_N$ for the
projection space in \eqref{def}. The results of calculation of
$\spec_2(H,\mathbf{M}_N\times\mathbf{L}_N)$ are shown in
Figures~\ref{fig:spec2_H_detail} and~\ref{fig:spec2_H}. The first
figure (with different scales along the real and imaginary axes)
demonstrates once more the advantage of this method in eliminating
the spectral pollution. Note also the series of points of
$\spec_2(H,\mathbf{M}_N\times\mathbf{L}_N)$ along the line
$\re\lambda=1$ --- this indicates the presence of the essential
spectrum of $H$ at $\lambda = 1$. The second figure shows all
computed points of the second order relative spectra for $N=25, 50$, 
and $100$.

The points of the second order relative
spectra give good approximations for the low eigenvalues of $H$. We
postpone the discussion of the rates of convergence in various
cases until Section 7. We note only that the number of well approximated 
eigenvalues increases with $N$, see Figure~\ref{fig:spec2_H}.

Examples I and II are discussed further in Sections 6 and 7; we
proceed now to the proofs of the main theoretical results stated
in Section 2.

\clearpage

\section{Proofs of the main results}

\begin{proof}[Proof of Theorem \ref{poll}]
It is clear that there exists
$\lambda_- \in \spech_{\text{ess}}(A)$
which lies strictly to the left of $\lambda$. Consider
first the case $\lambda_- \in \spec_{\text{ess}}(A)$.
If $\lambda_-$ is an isolated point of
$\spec(A)$, then it is an eigenvalue of infinite
multiplicity and we define
$\mathcal{H}^-_j := \mbox{span}\left(\{e_j\}\right)$,
where $e_j$, $j \in \mathbb{N}$ are mutually orthogonal eigenvectors
of $A$ corresponding to $\lambda_-$.
If $\lambda_-$ is not an isolated point, there exist
pair-wise disjoint compact sets
$W^-_j \subset \mathbb{R}$ lying  strictly to the left of
$\lambda$ and such that
$W^-_j\cap\spec(A) \not= \varnothing$, $j \in \mathbb{N}$.
In this case we choose $\mathcal{H}^-_j = E(W^-_j)\mathcal{H}$,
where $E(\cdot)$ are the spectral projections corresponding to $A$.
Finally,
if no point of $\spec_{\text{ess}}(A)$
lies strictly to the left of $\lambda$, then
$\lambda_- = -\infty$ and there exist eigenvalues
$\lambda_j < \lambda$ of $A$ converging to
$-\infty$. Here we set
$\mathcal{H}^-_j := E(\{\lambda_j\})$.
In all the above cases
$\mathcal{H}^-_j \subset E((-\infty, \lambda))$,
$j \in \mathbb{N}$.

Similarly one can define subspaces
$\mathcal{H}^+_j \subset E((\lambda, +\infty))$,
$j \in \mathbb{N}$. It is clear that $\mathcal{H}^\pm_j$
are all mutually orthogonal,
$\mathcal{H}^\pm_j \subset \dom(A)$,
$A\mathcal{H}^\pm_j \subset \mathcal{H}^\pm_j$,
and there exist $r_k > 0$ such that $r_k \nearrow +\infty$
and
$\mathcal{H}^\pm_j \subset E([-r_k, r_k])\mathcal{H}$,
$j \le k$.
Let $\mathcal{H}_0$ be the closed linear subspace
spanned by $\mathcal{H}^\pm_j$, $j \in \mathbb{N}$
and $\mathcal{H}_0^\perp$ be its orthogonal complement.

Let $v^\pm_k \in \mathcal{H}^\pm_k$, $\|v^\pm_k\| = 1$
and $\mu^\pm_k := (Av^\pm_k, v^\pm_k)$.
Then it follows from the spectral theorem that
$\mu^-_k < \lambda < \mu^+_k$, $k \in \mathbb{N}$.
Let $\mathcal{L}_k$ be the closed linear subspace
spanned by $\mathcal{H}^\pm_j$, $j < k$, \
$E([-r_k, r_k])\mathcal{H}_0^\perp$, and the vector
$$
u_k :=
\sqrt{\frac{\mu^+_k - \lambda}{\mu^+_k - \mu^-_k}}\
v^-_k +
\sqrt{\frac{\lambda - \mu^-_k}{\mu^+_k - \mu^-_k}}\
v^+_k .
$$
It is easily seen that the sequence
$(\mathcal{L}_k)_{k \in \mathbb{N}}$ is increasing
and belongs to $\Lambda(A)$. Since
$v^\pm_k, Av^\pm_k \in \mathcal{H}^\pm_k \subset
E([-r_k, r_k])\mathcal{H}$,
the orthogonality properties of
$\mathcal{H}^\pm_j$ imply
\begin{eqnarray}
& & \|u_k\| = 1 , \nonumber \\
& & ((A - \lambda I)u_k, u_k) =
\frac{\mu^+_k - \lambda}{\mu^+_k - \mu^-_k} (Av^-_k, v^-_k) +
\frac{\lambda - \mu^-_k}{\mu^+_k - \mu^-_k} (Av^+_k, v^+_k) -
\nonumber \\
& & \lambda \|u_k\|^2 =
\frac{\mu^+_k - \lambda}{\mu^+_k - \mu^-_k} \mu^-_k +
\frac{\lambda - \mu^-_k}{\mu^+_k - \mu^-_k} \mu^+_k -
\lambda = 0 , \nonumber \\
& & ((A - \lambda I)u_k, v) = 0 , \ \
\forall v \in E([-r_k, r_k])\mathcal{H}_0^\perp\cup
\bigcup_{j < k} (\mathcal{H}^-_j\cup\mathcal{H}^+_j) .
\nonumber
\end{eqnarray}
Hence $((A - \lambda I)u_k, u) = 0$, $\forall u \in
\mathcal{L}_k$, i.e. $P_k(A - \lambda I)u_k = 0$,
i.e. $P_kAu_k = \lambda u_k$, i.e. $\lambda$
is an eigenvalue of $P_kA : \mathcal{L}_k \to \mathcal{L}_k$,
$\forall k \in \mathbb{N}$.
\end{proof}

\begin{rem} If $\mathcal{H}$ is separable, the subspaces
$\mathcal{L}_k$ in Theorem \ref{poll} can be made finite
dimensional. Indeed, let $e^\pm_{j, m}$, $m = 1, \dots$
be an orthonormal basis of $\mathcal{H}^\pm_j$ and
$e^0_m$, $m = 1, \dots$ be an orthonormal basis of
$\mathcal{H}_0^\perp$. These bases exist because all subspaces
of $\mathcal{H}$ are separable.
Let $\mathcal{L}_k$ be the closed linear subspace
spanned by $e^\pm_{j, m}$ and $E([-r_k, r_k])e^0_m$ with  $j, m < k$,
and the vector
$$
u_k :=
\sqrt{\frac{\mu^+_k - \lambda}{\mu^+_k - \mu^-_k}}\
e^-_{k, 1} +
\sqrt{\frac{\lambda - \mu^-_k}{\mu^+_k - \mu^-_k}}\
e^+_{k, 1} .
$$
Then it is clear that $\mathcal{L}_k$ is finite dimensional, and the
rest of the proof of  Theorem \ref{poll} remains unchanged.
\end{rem}

\begin{proof}[Proof of Theorem \ref{unb}]
It is sufficient to prove that for any finite dimensional
$\mathcal{L} \subset \dom(A)$ and any
$\varepsilon, R > 0$ there exists  $\mathcal{L}' \subset \dom(A)$
such that
$$
\|P- P'\| < \varepsilon \ \ \ \mbox{and} \ \ \
\spec(A, \mathcal{L}') \subset (R, +\infty) ,
$$
where $P : \mathcal{H} \to \mathcal{L}$ and
$P' : \mathcal{H} \to \mathcal{L}'$ are the corresponding
orthogonal projections.

Let $m := \dim \mathcal{L}$ and $\varphi_1, \dots, \varphi_m$
be an orthonormal basis of $\mathcal{L}$. Since $A$ is not bounded
above, for any $\tau > 0$ there exist mutually
orthogonal vectors $\psi_1, \dots, \psi_m \in
\dom(A)\cap E((\tau, +\infty))\mathcal{H}$ which
we normalize by
$\|\psi_l\| = \delta$, $l = 1, \dots, m$,
where $\delta \in (0, 1)$. The numbers $\tau$ and
$\delta$ will be chosen later.

It is easy to see that the vectors
$\varphi_1 + \psi_1, \dots, \varphi_m + \psi_m$ are linearly independent.
Indeed, for any $c_1, \dots, c_m \in \mathbb{C}$ Parseval's identity
implies
\begin{equation}\label{lower}
\left\|\sum_{l = 1}^m c_l (\varphi_l + \psi_l)\right\| \ge
\left\|\sum_{l = 1}^m c_l \varphi_l\right\| -
\left\|\sum_{l = 1}^m c_l \psi_l\right\| =
(1 - \delta)\left(\sum_{l = 1}^m |c_l|^2\right)^{1/2} .
\end{equation}

Similarly
\begin{equation}\label{upper}
\left\|\sum_{l = 1}^m c_l (\varphi_l + \psi_l)\right\| \le
(1 + \delta)\left(\sum_{l = 1}^m |c_l|^2\right)^{1/2} .
\end{equation}

Let $\mathcal{L}'$ be the subspace spanned by
$\varphi_1 + \psi_1, \dots, \varphi_m + \psi_m$ and consider the
operators $J : \mathcal{L} \to \mathcal{L}'$ and
$J' : \mathcal{L}' \to \mathcal{L}$ defined by
$$
J\left(\sum_{l = 1}^m c_l \varphi_l\right) =
\sum_{l = 1}^m c_l (\varphi_l + \psi_l) , \ \
J'\left(\sum_{l = 1}^m c_l (\varphi_l + \psi_l)\right) =
\sum_{l = 1}^m c_l \varphi_l .
$$
It is clear that $J'J = I_{\mathcal{L}}$, $JJ' = I_{\mathcal{L}'}$
and
$$
(I - J)\left(\sum_{l = 1}^m c_l \varphi_l\right) =
-\sum_{l = 1}^m c_l \psi_l , \ \
(I - J')\left(\sum_{l = 1}^m c_l (\varphi_l + \psi_l)\right) =
\sum_{l = 1}^m c_l \psi_l .
$$
It follows from Parseval's identity and \eqref{lower},
\eqref{upper} that
$$
\|J\| \le 1 + \delta , \ \  \|J'\| \le \frac1{1 - \delta} , \ \
\|(I - J)P\| \le \delta , \ \  \|(I - J')P'\| \le \frac\delta{1 - \delta}\, .
$$
Further,
$$
\|P - P'\| \le \|P - P'P\| + \|P'P - P'\| = \|P - P'P\| + \|PP' - P'\| ,
$$
since $(P'P - P')^* = PP' - P'$. Further, for any $x \in \mathcal{H}$
the closest point of $\mathcal{L}'$ is $P'x$. Hence
$\|P - P'P\| \le \|P - JP\|$. Similarly $\|P' - PP'\| \le
\|P' - J'P'\|$. Consequently
$$
\|P - P'\| \le \|(I - J)P\| + \|(I - J')P'\| \le \delta +
\frac\delta{1 - \delta} = \frac{\delta(2 - \delta)}{1 - \delta} .
$$
Fix  $\delta \in (0, 1)$ such that
$\delta(2 - \delta)/(1 - \delta) < \varepsilon$. Then
$\|P- P'\| < \varepsilon$.

Since $\mathcal{L}$ is a finite dimensional subspace of
$\dom(A)$, there exists $M > 0$ such that
$\|Au\| \le M\|u\|$, $\forall u \in \mathcal{L}$.
Take an arbitrary
$$
v = \sum_{l = 1}^m c_l (\varphi_l + \psi_l) \in
\mathcal{L}'
$$
and set
$$
\varphi := \sum_{l = 1}^m c_l \varphi_l \in \mathcal{L} , \ \
\psi := \sum_{l = 1}^m c_l \psi_l \in
\dom(A)\cap E((\tau, +\infty))\mathcal{H} .
$$
Then
\begin{eqnarray}
& & (Av, v) = (A\psi, \psi) + 2\re (A\varphi, \psi) +
(A\varphi, \varphi) \ge \nonumber \\
& & \tau\|\psi\|^2 - 2M\|\varphi\| \|\psi\| -
M\|\varphi\|^2 =
(\tau\delta^2 - 2M\delta - M) \sum_{l = 1}^m |c_l|^2 .
\nonumber
\end{eqnarray}
Choose $\tau$ such that
$(\tau\delta^2 - 2M\delta - M)/(1 + \delta)^2 \ge R + 1$. Then
$$
(Av, v) \ge \frac{\tau\delta^2 - 2M\delta - M}{(1 + \delta)^2}
\|v\|^2 \ge (R + 1)\|v\|^2
$$
(see \eqref{upper}). Hence
$\spec(A, \mathcal{L}') \subset [R + 1, +\infty)$.
\end{proof}

Suppose $A$ is a bounded below
self-adjoint operator. Then everything said above generalizes from
$\Lambda(A)$ to $\Lambda(Q_A)$ (see Section 2 for the definition)
with almost no changes. One only needs to modify
the concluding estimates of the proof of Theorem \ref{unb}
in the following way. Adding if necessary a sufficiently
large positive number to the bounded below operator $A$,
one can assume without loss of generality that
$A$ is positive definite. Then
\begin{eqnarray}
& & Q_A(v, v) = Q_A(\psi, \psi) + 2\re \, Q_A(\varphi, \psi) +
Q_A(\varphi, \varphi) \ge \nonumber \\
& & Q_A(\psi, \psi) -
2Q_A(\psi, \psi)^{1/2} Q_A(\varphi, \varphi)^{1/2}  \ge \nonumber \\
& & \frac12 Q_A(\psi, \psi) - 2Q_A(\varphi, \varphi) \ge
\frac{\tau}2 \|\psi\|^2 -
2 M\|\varphi\|^2 = \nonumber \\
& & \frac{\tau\delta^2 - 4M}2 \sum_{l = 1}^m |c_l|^2 .
\nonumber
\end{eqnarray}
Choose $\tau$ such that
$(\tau\delta^2 - 4M)/2(1 + \delta)^2 \ge R + 1$. Then
$$
Q_A(v, v) \ge \frac{\tau\delta^2 - 4M}{2(1 + \delta)^2}
\|v\|^2 \ge (R + 1)\|v\|^2 .
$$
Hence
$\spec(A, \mathcal{L}') \subset [R + 1, +\infty)$.

Now we proceed to the proof of Theorem~\ref{mth}. For any $c_1,
c_2 \in \mathbb{C}$ we denote by $\mathbb{D}(c_1, c_2)$ the open
disk with diameter $(c_1, c_2)$.

\begin{lemma}\label{main}
Let $A$ be a self-adjoint operator on a Hilbert space $\mathcal{H}$
and let $\mathcal{L}$ be a finite dimensional subspace of
$\dom(A)$.
Suppose $a, b \in \mathbb{R}$ are such that
$$
(a, b) \cap \spec(A) = \varnothing .
$$
Then
$$
\mathbb{D}(a, b) \cap \spec_2(A, \mathcal{L})
= \varnothing .
$$
\end{lemma}

\begin{proof}[Proof of Lemma~\ref{main}]
Take an arbitrary $z \in \mathbb{D}(a, b)$.
The interval $(a, b)$ is seen from $z$ at an angle
$\pi/2 + \varepsilon$, where
$\varepsilon \in (0, \pi/2]$.
It is clear that the set
$$
F = \left\{(\lambda - z)^2 : \
\lambda \in \spec(A)\right\}
$$
lies in a sector with vertex at 0 and angle $\pi - 2\varepsilon$.
The distance from 0 to $F$ is obviously positive. Therefore there
exist $\vartheta \in (-\pi, \pi]$ and $c > 0$ such that
$$
\re \left(e^{i\vartheta}(\lambda - z)^2\right) \ge c, \
\forall \lambda \in \spec(A) .
$$
Now applying the spectral theorem for $A$ one obtains
\begin{eqnarray}
& & \re \Big(e^{i\vartheta}
\left((A - z I)u, (A - \overline{z} I)u\right)\Big)
= \nonumber \\
& & \re \Big(e^{i\vartheta}
\left(\|Au\|^2 - 2z (Au, u) + z^2\|u\|^2\right)\Big) =
\nonumber \\
& & \re \left(e^{i\vartheta}
\int_{\mbox{{\scriptsize Spec}}(A)}
(\lambda^2 -2\lambda z + z^2) d(E(\lambda)u, u)\right) =
\nonumber \\
& & \int_{\mbox{{\scriptsize Spec}}(A)}
\re \left(e^{i\vartheta}(\lambda - z)^2\right)
d(E(\lambda)u, u) \ge
c \int_{\mbox{{\scriptsize Spec}}(A)}  d(E(\lambda)u, u) =
\nonumber \\
& & c \|u\|^2 , \ \ \forall u \in \mathcal{L} .
\nonumber
\end{eqnarray}
Hence $u \in \mathcal{L}\setminus\{0\}$ cannot satisfy
\eqref{def}, i.e.
$z \not\in \spec_2(A, \mathcal{L})$.
\end{proof}

\begin{proof}[Proof of Theorem \ref{mth}]
Suppose the intersection in the statement of the Theorem is empty.
Then there exists $\varepsilon > 0$ such that
$$
\spec(A) \cap
[\re z - |\im z| - \varepsilon,
\re z + |\im z| + \varepsilon]
= \varnothing .
$$
Hence $\mathbb{D}(\re z - |\im z| - \varepsilon,
\re z + |\im z| + \varepsilon)$
does not intersect $\spec_2(A, \mathcal{L})$
according to Lemma \ref{main}. On the other hand, $z$ belongs
to both of these sets. The obtained contradiction proves the
theorem.
\end{proof}

\section{More on Example I}

Here we describe what happens if the function $a$ of Section 3 is
slightly modified. Since none of the changes we consider below
affects the second order relative spectra in any serious way, we
concentrate on the standard (first order) projection method.

We have performed calculations, similar to those presented in
Section 3, for the operator $A$ corresponding to the function
which is obtained from $a$ by multiplying the cosines  by 0.01,
0.001, or zero. The only noticeable difference with the results in
Section 3 is that this time there are more spurious eigenvalues in
the gap of the essential spectrum and additional spurious
eigenvalues appear, as $N$ increases, quicker than in the original
example (cf. Table~\ref{tab:spec_A} and Remark \ref{slow}).

The situation becomes more interesting if one moves the
discontinuity of $a$ from 0 to a point $y \not= 0$:
\begin{equation}\label{newa}
a(x) = a_y(x) =
\begin{cases}
\displaystyle -\frac{3}{2}+\frac{1}{2}\cos\sqrt{5}\,x\,,
&\qquad\text{for}\quad -\pi\le x< y\,,\\
\displaystyle
2+\cos\sqrt{2}\,x\,, &\qquad\text{for}\quad y\le x<\pi\,.
\end{cases}
\end{equation}
The corresponding operator of multiplication by $a_y$ will be
denoted $A_y$.

\begin{figure}[thb!]
\begin{center}
\resizebox{0.9\textwidth}{!}{\includegraphics*{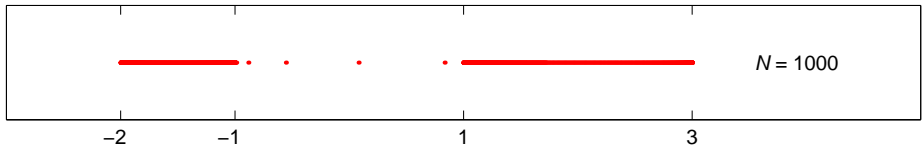}}
\end{center}
\caption{\label{fig:spec_A_L_2n}\small $\spec(A,\mathcal{L}_{1000})$.}
\end{figure}

To put the results into the right perspective, we need to return
to the original example $A=A_0$ first. We know from Section 3 that
$\spec(A, \mathcal{L}_N)$ show remarkable stability as $N = 2n +
1$ increases. The first surprise happens when one looks at
$\spec(A, \mathcal{L}_N)$ for $N = 2n$, where $\mathcal{L}_N :=
\text{span}\,(\{e_k\}_{k = -n}^{n - 1})$. These also seem to be
converging, but the spurious eigenvalues differ considerably from
those corresponding to $N = 2n + 1$ (see
Figure~\ref{fig:spec_A_L_2n} and Table~\ref{tab:spec_A_2n}  and
compare them with Figure~\ref{fig:spec_A_L_N} and
Table~\ref{tab:spec_A}).

\begin{table}[thb!]
\begin{center}
\begin{tabular}[t]{|c|*{4}{r}|}
\hline $N$ & 100     & 400     & 700     & 1000    \\ \hline
           &         & -0.9988 & -0.9935 & -0.9882 \\
           & -0.9910 & -0.9331 & -0.9033 & -0.8832 \\
$\spec(A,\mathcal{L}_N)\cap(-1,1)$ &  -0.7721&  -0.6358 & -0.5853 & -0.5546 \\
           & -0.0661 &  0.0317 &  0.0636 &  0.0821 \\
           &  0.8883 &  0.8539 &  0.8418 &  0.8345 \\ \hline
\end{tabular}
\end{center}
\caption{\label{tab:spec_A_2n}\small Eigenvalues lying in the gap
of the essential spectrum of $A$ for some even values of $N$.}
\end{table}

One might argue that this difference between the two cases occurs
due to the nonsymmetric truncation in the case $N = 2n$ and that
taking only odd values of $N$ is more natural. We prefer, however,
to interpret this fact as {\em $2$-periodicity in $N$ of the
points of $\spec(A, \mathcal{L}_N)$ lying in the gap $(-1,1)$
of the essential spectrum of $A$, i.e., of the spurious eigenvalues of $A$}.
Such an interpretation is supported by
Table~\ref{tab:spec_A_manyN}, which shows the spurious eigenvalues
for some consecutive values of $N$.

\begin{table}[thb!]
\begin{center}
{\scriptsize
\begin{tabular}[t]{|c|*{8}{r}|}
\hline $N$ & 301 & 302 & 303 & 304 & 305 & 306 & 307 & 308\\\hline
 &  -0.9899 & -0.9470 & -0.9898 & -0.9467 & -0.9896 & -0.9464 & -0.9895 & -0.9461\\
 &  -0.8468 & -0.6622 & -0.8463 & -0.6616 & -0.8458 & -0.6610 & -0.8453 & -0.6603\\
\raisebox{1.5ex}[0pt]{$\spec(A,\mathcal{L}_N)\cap(-1,1)$} &
    -0.3740 &  0.0143 & -0.3734 &  0.0147 & -0.3728 &  0.0151 & -0.3721 &  0.0155\\
 &   0.4554 &  0.8604 &  0.4555 &  0.8603 &  0.4556 &  0.8601 &  0.4557 &  0.8600\\\hline
\end{tabular}
}
\end{center}
\caption{\label{tab:spec_A_manyN}\small Eigenvalues lying in the gap
of the essential spectrum of $A$ for some consecutive values of $N$. 
2-periodicity in $N$ is clearly visible.}
\end{table}

We show below that in a more general case of the operator $A_y$
the ``period" depends on the location of the point $y$ of the
discontinuity of $a_y$ (see \eqref{newa}). This period is equal to
$2$ in the original example only because we take $y = 0$ there.

\begin{figure}[thb!]
\begin{center}
\resizebox{0.9\textwidth}{!}{\includegraphics*{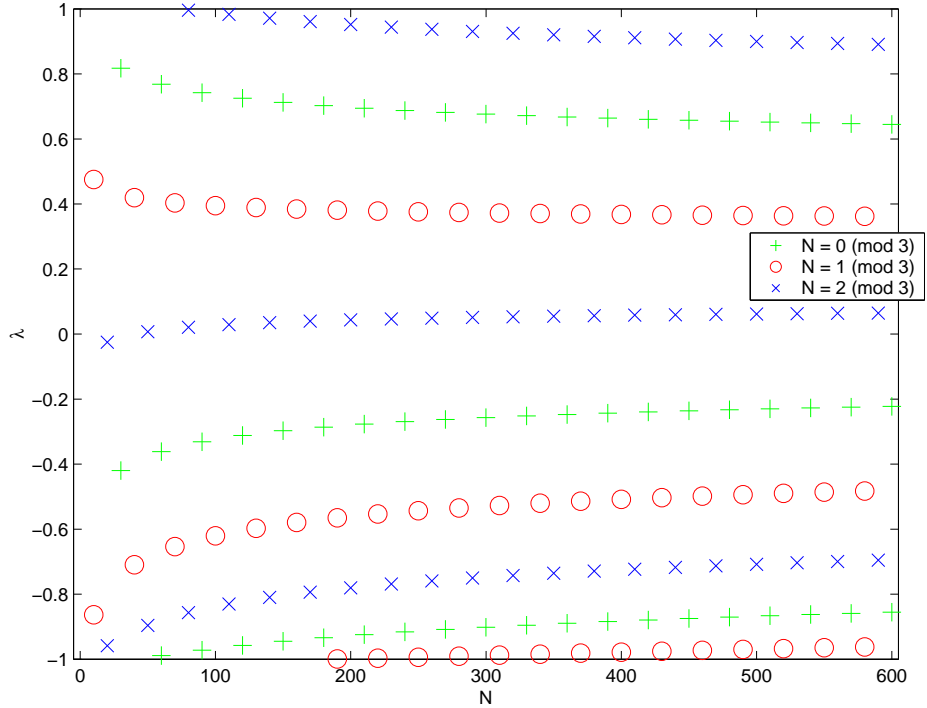}}
\end{center}
\caption{\label{fig:spec_A_z_versus_N}\small  Points of
$\spec(A_y,\mathcal{L}_N)$, $y=-\pi/3$, in the interval $(-1,1)$
as functions of $N$.}
\end{figure}

Let $y = -\pi/3$. Figure \ref{fig:spec_A_z_versus_N} shows the
computed (spurious) eigenvalues of $A_y = a_yI$ in the interval
$(-1,1)$ as functions of $N$. Remarkably, the eigenvalues vary
quite slowly with $N$ as long as $N=\text{const}\pmod{3}$. To
emphasize this ``3--periodicity" in $N$, we use different symbols
for points with different values of $N \pmod{3}$.

\begin{figure}[thb!]
\begin{center}
\resizebox{0.9\textwidth}{!}{\includegraphics*{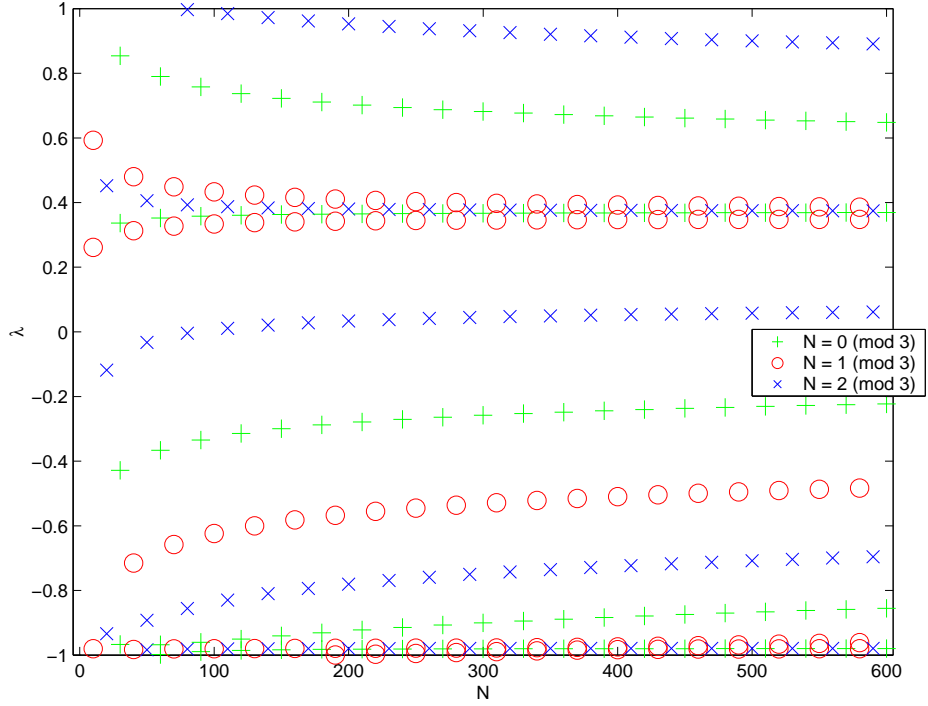}}
\end{center}
\caption{\label{fig:spec_B_z_versus_N}\small Points of
$\spec(B_y,\mathcal{L}_N)$, $y=-\pi/3$, in the interval $(-1,1)$
as functions of $N$. Note the presence of both ``genuine" and
``spurious" eigenvalues.}
\end{figure}

A similar effect occurs for the operator $B_y$, which is obtained
from the operator $B$ in \eqref{eq:B} by replacing $a$ with $a_y$
(see Figure \ref{fig:spec_B_z_versus_N}). However, curves
corresponding to the eigenvalues in the interval $(-1,1)$, which
were confirmed by the second order relative spectra as the correct
ones, do not depend on the value of $N  \pmod{3}$. They also seem
to converge to constants faster than those depending on $N\pmod{3}$.

\begin{figure}[thb!]
\begin{center}
\resizebox{0.9\textwidth}{!}{\includegraphics*{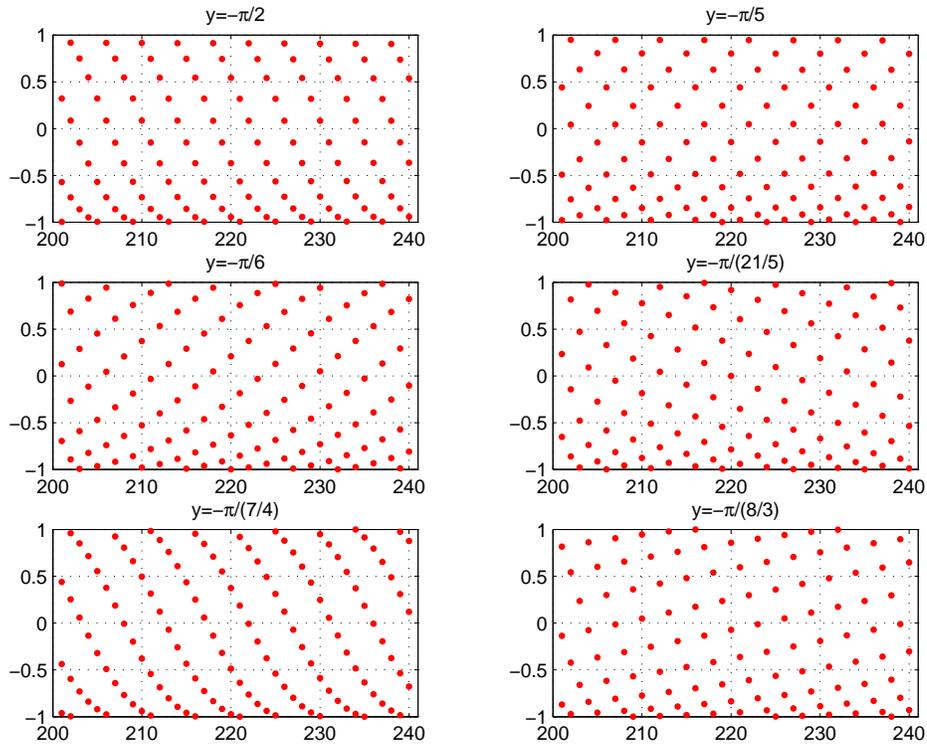}}
\end{center}
\caption{\label{fig:all_z}\small Points of
$\spec(A_y,\mathcal{L}_N)$ in the interval $(-1,1)$ as functions
of $N$, plotted for various values of $y$.}
\end{figure}

We have performed similar calculations with $y = -\pi/(p/q) = -\pi
q/p$ for the following values $p/q =$ 2, 4, 5, 6, 7, 21, 21/10,
21/5, 7/3, 8/3, and 5/2 (see Figure~\ref{fig:all_z} for the results for 
some of the values). In all these cases the spurious eigenvalues of $A_y$
showed a very clear pattern of ``periodic" behaviour with the
minimal ``period"
$$
\omega(q/p) :=
\begin{cases}
2, &\qquad\text{if}\quad q = 0,\\
p, &\qquad\text{if}\quad p \ \text{and} \ q \ \text{are both odd,}\\
2p, &\qquad\text{if}\quad q \not= 0 \ \text{and either} \
p \ \text{or} \ q \ \text{is even.}
\end{cases}
$$
Such ``periodic" behaviour can be observed in Figure~\ref{fig:all_z}. The
appearance of $\omega(q/p)$ here is probably not very surprising,
as it is, for any simple fraction $q/p$, the smallest natural
number such that $\omega\times\,$(length of $[-\pi, -\pi q/p]$)
and $\omega\times\,$(length of $[-\pi q/p, \pi]$) are integer
multiples of $2\pi$.

\section{More on Example II}

As we have shown in Section 4, a FEM calculation of the spectrum
of operator $H$ using the projection space
$\mathbf{L}_N\times\mathbf{L}_N$ leads to heavy spectral pollution
in the interval $(1,2)$. At the same time, the second order
spectra relative to the projection space
$\mathbf{M}_N\times\mathbf{L}_N$ approximate the spectrum of $H$
very well.

\begin{figure}[thb!]
\begin{center}
\resizebox{0.9\textwidth}{!}{\includegraphics*{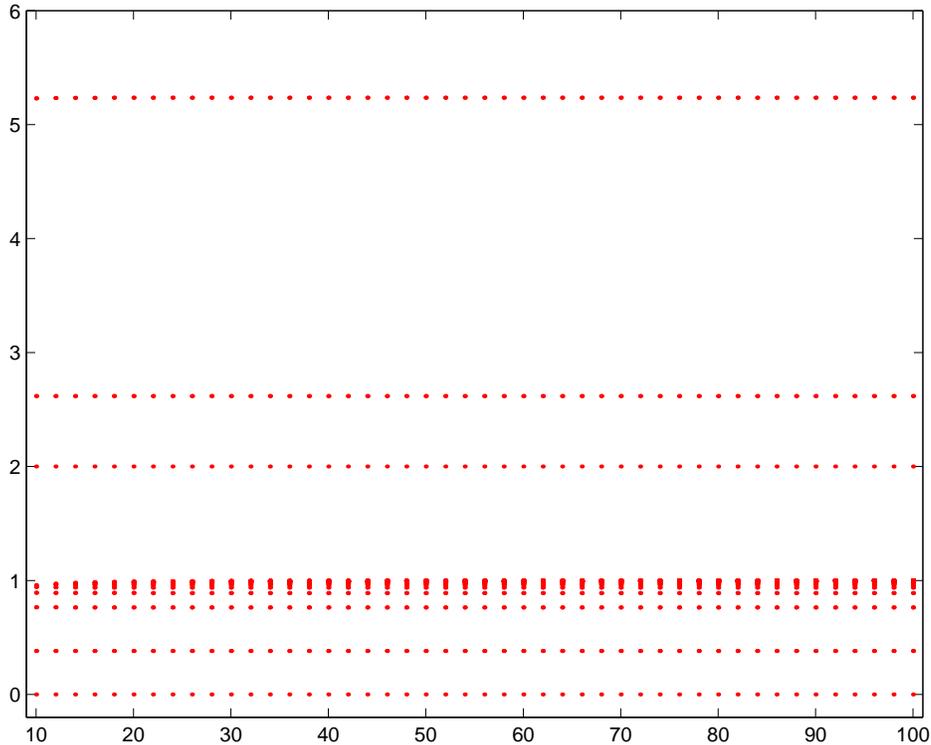}}
\end{center}
\caption{\label{fig:spec_H_ML}\small FEM calculated
$\spec(H,\mathbf{M}_N\times\mathbf{L}_N)$ plotted against $N$.
There are no spurious eigenvalues!}
\end{figure}

A natural question is whether such a discrepancy may be mostly due
to the {\em different} projection spaces being used in the
calculations of the first order and second order spectra. To check
this, we compute the spectra
$\spec(H,\mathbf{M}_N\times\mathbf{L}_N)$ using the same
projection spaces as we used in Section 4 for the second order relative
spectra. The results (restricted to eigenvalues lying in the
interval $[0,6]$) are shown in
Figure~\ref{fig:spec_H_ML}. This time, there  are no spurious
eigenvalues at all! In the particular case $N=100$, we have 100
eigenvalues in $[0,1]$, and 200 eigenvalues in $[2,\infty)$
(recall that $\dim(\mathbf{M}_N\times\mathbf{L}_N)=3N$). The
actual eigenvalues of $H$ are well approximated in this case as
well, see Remark~\ref{rem:accuracy} below.

This of course indicates that the projection spaces
$\mathbf{M}_N\times\mathbf{L}_N$ are much better suited for the
calculation of the spectrum of $H$ than the projection spaces
$\mathbf{L}_N\times\mathbf{L}_N$. Such a dependence upon the
choice of the spaces of finite elements is quite interesting,
especially because it does not occur for the
operator $-d^2/dx^2$ and because of the explicit relations which
exist between the eigenvalues and eigenfunctions of $H$ and those
of $-d^2/dx^2$, see \eqref{eq:eigsH} and \eqref{eigv}.

The spaces $\mathbf{M}_N\times\mathbf{L}_N$ are well suited for the
calculation of $\spec(H)$ because they are subspaces of the mixed
order Sobolev space $W_2^2([-\pi,\pi)])\times W_2^1([-\pi,\pi)])$ and
therefore reflect much better the {\em mixed orders of
differentiation} in the operator $H$, see also \cite{RSSV} and
references therein. We attempt to explain this phenomenon in a
more general setting below, see Remark~\ref{rem:mixed}.

\begin{figure}[thb!]
\begin{center}
\resizebox{0.9\textwidth}{!}{\includegraphics*{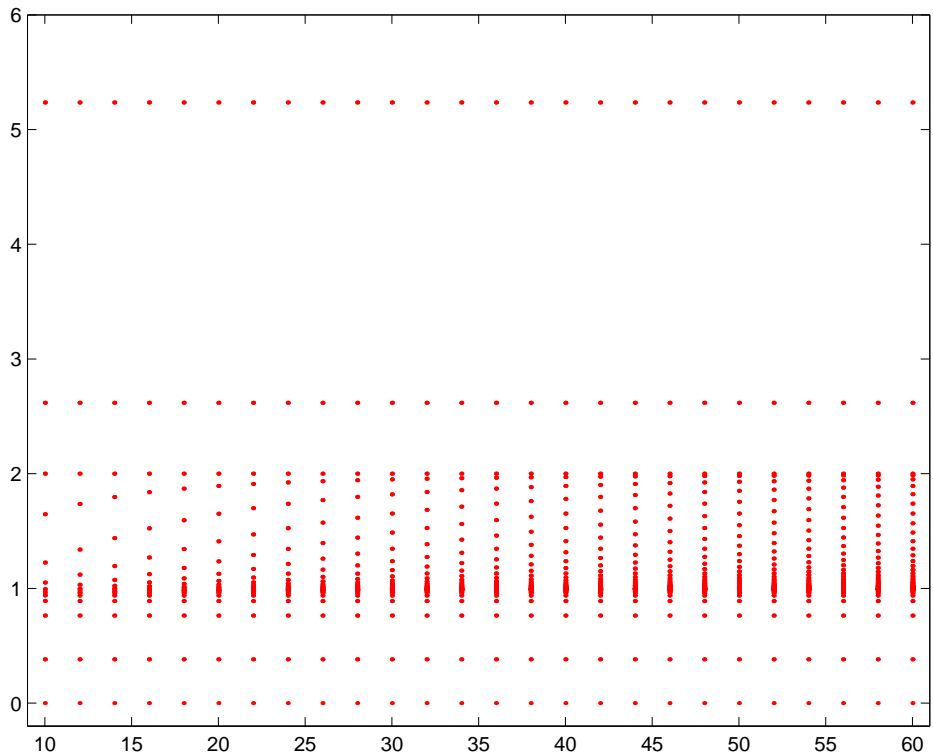}}
\end{center}
\caption{\label{fig:spec_H_MM}\small
$\spec(H,\mathbf{M}_N\times\mathbf{M}_N)$ plotted against $N$.
Again, all the eigenvalues in the interval $(1,2)$ are spurious!}
\end{figure}

One may be tempted to think that the success of the second order
method used in Section 4 is also only due to the correct choice of
the finite elements. In order to check this assumption let us
compare the results of calculation of the first order and second
order relative spectra using the projection spaces
$\mathbf{M}_N\times\mathbf{M}_N$. The results for 
$\spec(H,\mathbf{M}_N\times\mathbf{M}_N)$ are
shown in Figure~\ref{fig:spec_H_MM} --- again, there is spectral
pollution in the interval $(1,2)$.

\begin{figure}[thb!]
\begin{center}
\resizebox{0.9\textwidth}{!}{\includegraphics*{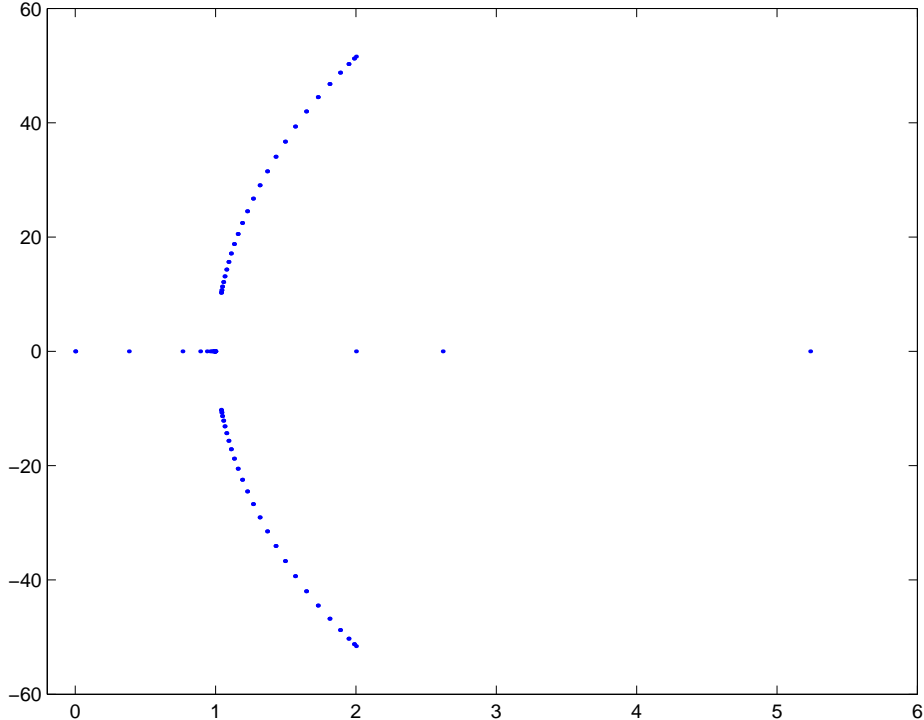}}
\end{center}
\caption{\label{fig:spec2_H_MM}\small
$\spec_2(H,\mathbf{M}_{50}\times\mathbf{M}_{50})$.}
\end{figure}

On the other hand, the calculations of the second order relative
spectra $\spec_2(H,\mathbf{M}_N\times\mathbf{M}_N)$ produce
pictures very similar to those in Figures~\ref{fig:spec2_H_detail}
and~\ref{fig:spec2_H}. The only noticeable difference is
the presence of a number of points with relatively large moduli of
imaginary parts ($\gtrsim 10$ for $N=50$) which seem to lie on a ``nice''
curve passing through the point $\lambda=1$ (see
Figure~\ref{fig:spec2_H_MM}; note the different scale along the
imaginary axis compared with either Figure~\ref{fig:spec2_H_detail} or
\ref{fig:spec2_H}).

Thus, for the spaces $\mathbf{M}_N\times\mathbf{M}_N$, the
projection method using second order relative spectra works quite
well, whereas the standard method again fails dramatically.

\begin{rem}\label{rem:accuracy} It should be noted that both the
first order and the second order relative spectra (outside
$(1,2)$) approximate $\spec(H)$ very well, with accuracy depending
on the finite elements chosen. The points of $\spec(H,
\mathbf{L}_N\times\mathbf{L}_N)$, $\spec(H,
\mathbf{M}_N\times\mathbf{L}_N)$, and $\spec(H,
\mathbf{M}_N\times\mathbf{M}_N)$ converge to the points
$\lambda_k^\pm\in\spec(H)$ with the rate of about $O(N^{-2})$,
$O(N^{-4})$, and $O(N^{-6})$, respectively. For the second and
third choice of finite elements, the same rates of convergence are
achieved by the real parts of the points of $\spec_2(H,
\mathbf{M}_N\times\mathbf{L}_N)$ and $\spec_2(H,
\mathbf{M}_N\times\mathbf{M}_N)$, whereas the imaginary parts
converge to zero with the rate of about $O(N^{-2})$ and
$O(N^{-3})$, respectively.
These convergence rates were obtained by analyzing the 
$\log$--$\log$ dependence of errors upon $N$ for several low 
eigenvalues taken from both $\lambda^+$ and $\lambda^-$ 
series.

\end{rem}

\begin{rem}\label{rem:mixed}
Here we try to explain  why the standard projection method of
approximating $\spec(H)$ does not lead to spectral pollution in
the case of the projection spaces
$\mathbf{M}_N\times\mathbf{L}_N$. The main reason is that it seems,
in this particular case, to ``mimic'' the way one finds $\spec(H)$
analytically. Indeed, let $\mathcal{L}_1$ be a space of finite
elements for $v$, $\mathcal{L}_2$ be a space of finite elements
for $u$, and let $\mathcal{L}'_2 := \{\psi' : \  \psi \in
\mathcal{L}_2\}$. Let $\varphi_j$, $j = 1, \dots, M$ be a basis of
$\mathcal{L}_1$ and $\psi_k$, $k = 1, \dots, N$ be a basis of
$\mathcal{L}_2$. Since the projection method does not depend on
the choice of bases, provided all the calculations are precise, we
can assume that the basis $\{\varphi_j\}$ is orthonormal. Then the
finite dimensional system for \eqref{np1} takes the form
$$
\left\{ \begin{array}{rl}
\sum_{k = 1}^N (\psi'_m, \psi'_k)u_k +
\sum_{j = 1}^M (\psi'_m, \varphi_j)v_j =
&  \lambda \sum_{k = 1}^N (\psi_m, \psi_k)u_k , \\
& m =  1, \dots, N ,\\
\sum_{k = 1}^N (\varphi_j, \psi'_k)u_k +
2v_j = &  \lambda v_j , \\ & j = 1, \dots, M .
\end{array} \right.
$$
This leads to the following system for the coefficients
$u_k$:
\begin{eqnarray}\label{np2}
\sum_{k = 1}^N (\psi'_m, \psi'_k)u_k +
\frac1{\lambda - 2}
\sum_{k = 1}^N \sum_{j = 1}^M (\psi'_m, \varphi_j)
(\varphi_j, \psi'_k)u_k = \nonumber \\
\lambda \sum_{k = 1}^N (\psi_m, \psi_k)u_k , \ \
m =  1, \dots, N .
\end{eqnarray}
Let $P_1$ be the orthogonal projection onto
$\mathcal{L}_1$. Then
$$
P_1\psi'_k = \sum_{j = 1}^M (\psi'_k, \varphi_j) \varphi_j
\ \ \mbox{and} \ \
(P_1\psi'_m, P_1\psi'_k) = \sum_{j = 1}^M (\psi'_m, \varphi_j)
(\varphi_j, \psi'_k) .
$$
Hence \eqref{np2} is equivalent to
\begin{eqnarray}\label{np3}
\sum_{k = 1}^N (\psi'_m, \psi'_k)u_k +
\frac1{\lambda - 2}
\sum_{k = 1}^N (P_1\psi'_m, P_1\psi'_k)u_k =
\lambda \sum_{k = 1}^N (\psi_m, \psi_k)u_k , \\
m =  1, \dots, N .
\nonumber
\end{eqnarray}

On the other hand, it is clear that \eqref{np1} implies
\begin{equation}\label{np4}
-u'' - \frac1{\lambda - 2} u''
 = \lambda u\,,
\end{equation}
which is equivalent to \eqref{eq:sixstar}.
So, the projection method applied to \eqref{np4} should not lead
to spectral pollution. This gives the following finite dimensional
system
\begin{eqnarray}\label{np5}
\sum_{k = 1}^N (\psi'_m, \psi'_k)u_k +
\frac1{\lambda - 2}
\sum_{k = 1}^N (\psi'_m, \psi'_k)u_k =
\lambda \sum_{k = 1}^N (\psi_m, \psi_k)u_k , \\
m =  1, \dots, N .
\nonumber
\end{eqnarray}

If $\mathcal{L}_1 \supseteq \mathcal{L}'_2$, then \eqref{np3} and
\eqref{np5} coincide and the projection method applied to
\eqref{np1} does not pollute. If $\mathcal{L}_1 \not\supseteq
\mathcal{L}'_2$, but $\mathcal{L}_1$ and $\mathcal{L}'_2$ are
``compatible", then it may well happen that \eqref{np3} mimics the
properties of \eqref{np5} at least partially, and the projection
method applied to \eqref{np1} does not lead to spectral pollution.
\end{rem}

\section{Concluding remarks}

The main message of the paper is that projection methods
using second order relative spectra may provide an
efficient tool for dealing with spectral pollution.
They are not much harder to implement than the standard
projection methods and are likely to give useful {\em a posteriori\/}
information (see Theorem \ref{mth}) about the location
of (at least a part of) the spectrum of a self-adjoint
operator. We would feel that the paper has achieved its goal
if it persuaded other people to try second order
relative spectra in their problems of interest. We hope that
the method we suggest here will work well in problems more
serious than the model examples considered above.

We would like to stress again that we do not claim that
second order relative spectra are the only way of addressing
spectral pollution. In our Example I, for instance,
one could argue that spectral pollution could have
been detected just by looking at pictures like Figure 3, without
any reference to the second order relative spectra. The argument
would be that a
genuine eigenfunction ``cannot" look like the one in Figure 3, and
hence the corresponding eigenvalue is spurious.
The ``periodicity" of spurious eigenvalues (see Section 6)
might also have been used to detect possible spectral
pollution. Indeed,
those approximate eigenvalues of $B_y$, which were
confirmed by the second order relative spectra as the
correct ones, do not exhibit ``periodicity" in $N$. They
are present for any sufficiently large value of $N$.
We have no
intention to disagree with these arguments.
Note however that both of these approaches
are based on the results of calculations using the standard
(first order) projection method, which does not see the
difference between $A$ and the Toeplitz operator $T_a$.
Consequently, they cannot distinguish between spectral
pollution for the former and non-pollution for the
latter. Second order relative spectra, on the other hand,
are free from this problem. Our
point is that it is very easy to fail to notice
spectral pollution unless one is very careful, and
that second order relative spectra seem to be very
efficient ``lie detectors".

Another perfectly legitimate argument is that one should
try to avoid spectral pollution by making a sensible choice of
the subspaces $\mathcal{L}_k$, rather than allow it to
happen and then try to detect it. The standard projection
method using $\mathbf{M}_N\times\mathbf{L}_N$, for
instance, works very well in our Example II (see
Section 7). One should, however, bear in mind that choosing
``right" subspaces when applying the standard projection
method in more difficult situations may prove non-trivial
and require a good understanding of the problem at hand.
Second order relative spectra, on the other hand, seem
to provide a more robust projection method which works well
even in the case of ``wrong" subspaces, when the standard one
fails. The case of $\mathbf{M}_N\times\mathbf{M}_N$ considered
in Section 7 is an example of such situation.

Finally, we refer the reader to the companion paper by E.~B.~Davies and 
M.~Plum~\cite{DP}, where a different method of treating spectral pollution 
is developed.

\section*{Acknowledgements} 
\addcontentsline{toc}{section}{Acknowledgements}
We would like to thank Brian Davies for
interest in our work and  very helpful discussions.
This work was partially supported by a grant from the EPSRC 
Spectral Theory Network.

\end{document}